\documentclass[12pt]{amsart}

\newcommand {\supplus}{\mathop{{\supset}\llap{\raise 
0.5pt\hbox{\normalfont\small+}\hskip 0.5pt}}} 

\newcommand {\subplus}{\mathop{{\subset}\llap{\raise 
0.5pt\hbox{\normalfont\small+}\hskip 0.5pt}}}  

\newcommand {\Cee}    {{\mathbb  C}}

\newcommand {\Nee}    {{\mathbb  N}}

\newcommand {\Zee}    {{\mathbb  Z}}

\newcommand {\fa}     {{\mathfrak{a}}}

\newcommand {\fab}    {{\mathfrak{ab}}} 
\newcommand {\fag}    {{\mathfrak{ag}}}

\newcommand {\fb}     {{\mathfrak{b}}}

\newcommand {\fg}     {{\mathfrak{g}}}    %
\newcommand {\fgl}    {{\mathfrak{gl}}}  %
\newcommand {\fh}     {{\mathfrak{h}}}

\newcommand {\fo}     {{\mathfrak{o}}}
\newcommand {\fosp}   {{\mathfrak{osp}}}
\newcommand {\fp}    {{\mathfrak{p}}}   %
\newcommand {\fpe}    {{\mathfrak{pe}}}   %

\newcommand {\fpo}    {{\mathfrak{po}}}

\newcommand {\fpsl}   {{\mathfrak{psl}}}

\newcommand {\fq}     {{\mathfrak{q}}}

\newcommand {\fs}     {{\mathfrak{s}}}
\newcommand {\fS}     {{\mathfrak{S}}}

\newcommand {\fsl}    {{\mathfrak{sl}}}

\newcommand {\fspe}   {{\mathfrak{spe}}}

\newcommand {\fsq}    {{\mathfrak{sq}}}

\newcommand {\fsvect} {{\mathfrak{svect}}}

\newcommand {\fvect}  {{\mathfrak{vect}}}   %

\newcommand {\cal} {\mathcal}

\def \opname#1#2%
  {\expandafter\newcommand \csname #1\endcsname {{\mathop{#2}\nolimits}}}

\newcommand{\rmname}[1]
  {\expandafter\newcommand \csname #1\endcsname {{\operatorname{#1}}}}

\newcommand{\rmnameii}[2]
  {\expandafter\newcommand \csname #1\endcsname {{\operatorname{#2}}}}

\rmname{act}
\rmname{Ad}
\rmname{Add}
\rmname{ad}
\rmname{Aff}
\rmname{Alt}
\rmname{alt}
\rmname{Ann}
\rmname{antidiag}
\rmname{Ber}
\rmname{ber}
\rmname{Br}
\rmname{card}
\rmname{ch}
\rmname{Char}
\rmname{cem}
\rmname{cj}
\rmname{Cliff}
\rmname{cntr}
\rmname{codim}
\rmname{coind}
\rmname{const}
\rmname{col}
\rmname{cork}
\rmname{cpr}
\rmname{diag}
\rmnameii{Div}{div}
\rmname{Def}
\rmname{Der}
\rmname{Dim}
\rmname{End}
\rmname{Even}
\rmname{Ext}
\rmname{gr}
\rmname{Hom}
\rmname{HT}
\rmnameii{Ht}{ht}
\rmname{hwt}
\rmname{Id}
\rmname{id}
\rmname{ind}
\rmname{Ind}
\rmname{Inf}
\rmname{irr}
\rmname{Le}
\rmname{Lie}
\rmname{lwt}
\rmname{mult}
\rmname{Mor}
\rmname{nm}
\rmname{Ob}
\rmname{Odd}
\rmname{ord}
\rmname{Osc}
\rmname{per}
\rmname{Pic}
\rmname{pr}
\rmname{pro}
\rmname{Prime}
\rmname{Proj}
\rmname{prt}
\rmname{pt}
\rmname{Q}
\rmname{qet}
\rmname{qtr}
\rmname{rd}
\rmname{rk}
\rmname{row}
\rmname{Res}
\rmname{salt}
\rmname{Sch}
\rmname{SBr}
\rmname{scalar}
\rmname{sch}
\rmname{sh}
\rmname{Ser}
\rmname{sign}
\rmname{Smbl}
\rmname{spin}
\rmname{ssym}
\rmname{str}
\rmname{st}
\rmname{sgn}
\rmname{sq}
\rmname{symm}
\rmname{supp}
\rmname{Supp}
\rmname{St}
\rmname{Spec}
\rmname{Specm}
\rmname{Spm}
\rmname{tr}
\rmname{vpt}
\rmname{weyl}
\rmname{Weyl}
\rmname{Witt}

\opname{vvol}  {{v\hspace{-0.1ex}o\hspace{-0.02ex}l\/}}
\opname{pnt}  {\text{\normalfont pt}}
\opname{Span} {{Span}}
\opname{slim} {\overline{\lim}}
\opname{Vol}  {{V\hspace{-0.55ex}o\hspace{-0.02ex}l\/}}
\opname{QVol} {{Q\hspace{-0.3ex}V\hspace{-0.55ex}o\hspace{-0.02ex}l\/}}
\opname{PoVol}{{P\hspace{-0.35ex}o\hspace{-0.25ex}V\hspace{-0.55ex}o\hspace{-0.
02ex}l\/}}
\opname{BVol} {{B\hspace{-0.2ex}V\hspace{-0.55ex}o\hspace{-0.02ex}l\/}}
\opname{Par}  {{P\hspace{-0.3ex}a\hspace{-0.05ex}r\/}}

\rmname{Mat}
\rmname{Bil}
\rmname{Diff}
\rmname{Ker}
\rmname{Herm}
\rmname{Coker}
\rmname{Conn}
\rmname{Covect}
\rmname{Vect}
\rmname{Int}

\rmnameii {IM} {Im}
\rmnameii {RE} {Re}

\opname{Aut} {{A\hspace{-0.2ex}u\hspace{-0.1ex}t\/}}
\opname{GL} {{G\hspace{-0.3ex}L}}
\opname{SL} {{S\hspace{-0.3ex}L}}
\opname{Exp} {{E\hspace{-0.2ex}x\hspace{-0.1ex}p\/}}
\opname{GQ} {{G\hspace{-0.2ex}Q}}
\opname{OSp} {{O\hspace{-0.25ex}S\hspace{-0.15ex}p\/}}
\opname{Out} {{O\hspace{-0.25ex}u\hspace{-0.15ex}t\/}}
\opname{Spp} {{S\hspace{-0.2ex}p\/}}
\opname{SpO} {{S\hspace{-0.2ex}p\hspace{-0.02ex}O\/}}
\opname{Pe} {{P\hspace{-0.25ex}e\/}}
\opname{SPe} {{S\hspace{-0.25ex}P\hspace{-0.25ex}e\/}}
\opname{Spin} {{S\hspace{-0.25ex}p\hspace{-0.05ex}i\hspace{-0.1ex}n\/}}
\opname{Iso} {{I\hspace{-0.25ex}s\hspace{-0.1ex}o\/}}
\opname{SSPe} {{S\hspace{-0.25ex}S\hspace{-0.15ex}P\hspace{-0.25ex}e\/}}
\opname{PeU} {{P\hspace{-0.25ex}e\hspace{-0.1ex}U\/}}
\opname{QU} {{Q\hspace{-0.15ex}U\/}}
\opname{U} {{U\/}}

\opname{cGQ} {{\cal G \hspace{-0.2em} Q \/}}
\opname{cSL} {{\cal S \hspace{-0.2em} L \/}}
\opname{cGL} {{\cal G \hspace{-0.2em} L \/}}
\opname{cOSp} {{\cal O \hspace{-0.2em} S \hspace{-0.3em} \it p\/}}
\opname{cPe} {{\cal P \hspace{-1.5pt} \it e\/}}
\opname{cVect} {{\cal V \hspace{-1.5pt} \it e\hspace{-0.1ex}c\hspace{-0.1ex}t\/
}}
\opname{cVol} {{\cal V \hspace{-1.5pt} \it o\hspace{-0.1ex}l\/}}
\opname{cAut} {{\cal A \hspace{-0.2em} \it u\hspace{-0.1em}t\/}}
\opname{cCovect} {{\cal C \hspace{-1.5pt}
     \it o\hspace{-0.1ex}v\hspace{-0.1ex}e\hspace{-0.1ex}c\hspace{-0.1ex}t\/}}
\opname{CW} {{C\hspace{-0.15ex}W}}

\opname {Ab}   {{\sf Ab}}
\opname {Alg}   {{\sf Alg}}
\opname {ASch}  {{\sf Aff\;Sch}}
\opname {Funct}   {{\sf Funct}}
\opname {Gr}   {{\sf Gr}}
\opname {Grf}  {{{\sf Gr}_f}}
\opname {Mods}   {{\sf Mods}}
\opname {SMods}   {{\sf SMods}}
\opname {Rings}   {{\sf Rings}}
\opname {Salg}   {{\sf Salg}}
\opname {Sets} {{\sf Sets}}
\opname {SSMan} {{\sf SMan}}
\opname {Top}  {{\sf Top}}
\opname {Vebun}   {{\sf Vebun}}

\newcommand {\ev} {{\bar0}}
\newcommand {\od} {{\bar1}}
\newcommand {\eps} {\varepsilon}

\newcommand {\tto} {\longrightarrow}

\newcommand {\pderf}[2] {{\frac{\partial {#1}}{\partial {#2}}}}

\newcommand {\bcdot}   {\mathbin{\hbox{\raise.4ex\hbox{\bf.}}}} 

\newcommand {\secno} {}
\newcommand {\ssecfont} {\normalfont\bf}

\newtheorem{Theorem}{\secno Theorem}

\newtheorem{Corollary}[Theorem]{\secno Corollary}

\newtheorem{Problem}[Theorem]{\secno Problem}

\newenvironment {th*}[1]
    {\gdef\thname{#1} \begin{thn}}%
    {\end{thn}}
\newtheorem{thn}[Theorem] {\thname}

\theoremstyle{definition}

\newenvironment {ex*}[1]
    {\gdef\thname{#1} \begin{exn}}%
    {\end{exn}}
\newtheorem{exn}[Theorem]{\thname}

\theoremstyle{remark}

\newtheorem{Remark}[Theorem]{\secno Remark}

\newenvironment {rem*}[1]
    {\gdef\thname{#1} \begin{remn}}%
    {\end{remn}}
\newtheorem{remn}[Theorem]{\thname}

\newcommand {\ssec}{\subsection*}
\newcommand {\ssbegin}[2]
  {\def \secno {\gdef \secno {}{\ssecfont #1. }}%
   \begin{#2}}
\begin{document}

\title{The invariant polynomials on simple Lie superalgebras}

\author{Alexander Sergeev} 

\address{On leave of absence from the Balakovo Inst. of Technique of 
Technology and Control, Branch of Saratov State Technical University.  
Correspondence: c/o D. Leites, Dept.  of Math., Univ.  of Stockholm, 
Roslagsv.  101, Kr\"aftriket hus 6, S-106 91, Stockholm, Sweden}

\thanks{I am thankful to D.~Leites for help and support.}


\begin{abstract} Chevalley's theorem states that for any simple 
finite dimensional Lie algebra $\fg$ (1) the restriction homomorphism of 
the algebra of polynomials $S(\fg^*)\tto S(\fh^*)$ onto the Cartan 
subalgebra $\fh$ induces an isomorphism $S(\fg^*)^{\fg}\cong 
S(\fh^*)^{W}$, where $W$ is the Weyl group of $\fg$, (2) each 
$\fg$-invariant polynomial is a linear combination of the polynomials 
$\tr \rho(x)^k$, where $\rho$ is a finite dimensional representation 
of $\fg$.

None of these facts is necessarily true for simple Lie superalgebras.  We 
reformulate Chevalley's theorem as formula $(*)$ below to embrace Lie 
superalgebras.  Let $\fh$ be the split Cartan subalgebra of $\fg$; let 
$R=R_+\cup R_-$ be the set of nonzero roots of $\fg$, the union of 
positive and negative ones.  Set $\tilde R_+=\{\alpha \in R_+\mid 
-\alpha \in R_-\}$.  For each root $\alpha \in \tilde R_+$ denote by 
$\fg(\alpha)$ the Lie superlagebra generated by $\fh$ and the root 
superspaces $\fg_\alpha$ and $\fg_{-\alpha}$.  Let the image of 
$S(\fg(\alpha)^*)^{\fg(\alpha)}$ under the restriction homomorphism 
$S(\fg(\alpha)^*)\tto S(\fh^*)$ be denoted by $I^{\alpha}(\fh^*)$ and the 
image of $S(\fg^*)^{\fg}$ by $I(\fh^*)$.  Then
$$
I(\fh^*)=\mathop{\cap}\limits_{\alpha\in \tilde
R_+}I^{\alpha}(\fh^*).\eqno{(*)}
$$

Chevalley's theorem for anti-invariant polynomials is also given.
\end{abstract}

\maketitle

\section*{Introduction} The ground field is $\Cee$.  All Lie algebras 
and superalgebras and their representations are assumed to be of 
finite dimension unless stated to the contrary.

\ssec{0.1.  Lie algebras} The algebras of invariant polynomials on 
semisimple Lie algebras $\fg$ play an important role: in paricular, 
they allow one to describe the structure of the center $Z(\fg)$ of the 
universal envelopping algebra $U(\fg)$ of and deduce the Weyl 
character formula for the finite dimensional irreducible 
representations; the characteristic classes of vector bundles, analogs 
of Euler's equation of the solid body, etc.  are expressed in terms of 
the invariant polynomials.

If $\fg$ is semisimple, then it possesses an invariant nondegenerate 
symmetric bilinear form which induces an isomorphism $\fg\simeq 
\fg^{*}$.  This isomorphism commutes with the $\fg$-action and reduces 
the description of the center to the description of $\fg$-invariant 
polynomials on $\fg$, i.e., the elements of $S(\fg^{*})^\fg$.

The study of the algebras $S(\fg^{*})^\fg$ has the following  
important aspects:

1) If $\fg$ is semisimple, then having described $S(\fg^{*})^\fg$ we 
simultaneously describe $Z(\fg)$: these algebras are isomorphic.

2) Let $\fg$ be semisimple, $\fh$ its Cartan subalgebra, $W$ the Weyl 
group.  Then {\it Chevalley's theorem} states that the restriction 
homomorphism $S(\fg^*)\tto S(\fh^*)$ onto the Cartan subalgebra 
induces an isomorphism of the algebra of $\fg$-invariant polynomials 
on $\fg$ with the algebra of $W$-invariant polynomials on $\fh$:
$$
S(\fg^*)^{\fg}\cong S(\fh^*)^{W}.
$$
Moreover, each invariant polynomial of degree $k$ is a linear 
combination of the functions
$$
\tr \rho(x)^k,
$$
where $\rho$ runs over finite dimensional representations of $\fg$.

3) The algebras of $W$-invariant polynomials are important examples of 
of the algebras of invariant functions: for example, for $\fg=\fsl(n)$ 
this leads to the classical invariant theory of symmetric functions.

4) for a semi-simple Lie algebra there exists a homomorphism (called 
{\it Harish-Candra homomorphism})
$$
Z(\fg)\tto S(\fh)\simeq S(\fh^{*})\eqno{(0.1)}
$$
The image of the center under this homomorphism is the so-called 
algebra of {\it shifted} symmetric functions.  The study of functions 
of this type turned out to be important in the study of the {\it 
quantum immanents} (that constitute a distinguished basis of $Z(\fg)$).  
Its further generalization is connected with Macdonald polynomials.

\ssec{0.2.  Lie superalgebras} This paper is an attempt to realize 
part of the above program (1) -- (4) for simple Lie superalgebras.  
Let me point out the main distinctions:

i) not every simple Lie superalgebra possesses a nondegenerate 
symmetric invariant bilinear form. Even if it does possess one, the 
form might be odd, while for isomorphism $(0.1)$ we need an even 
form.

ii) The notion of the Weyl group becomes very involved in the 
supersetting.  In various instances it has to be described 
differently; moreover, even for the same problem the description of 
the analog differs drastically from one type of Lie superalgebras to 
another one.

iii) not every irreducible $\fg$-module is uniquely determined by its 
central character (this resembles the case of infinite dimensional 
modules over classical Lie algebras); hence, the central characters 
--- i.e., the invariant polynomials --- provide us with the character 
formula only for the generic --- {\it typical} --- representations;

iv) for several Lie superalgebras not every invariant polynomial can be
represented as a linear combination of supertraces of any finite
dimensional representations.
 
Recall that there are several types of simple Lie 
superalgebras with quite distinct properties:

a) possessing a Cartan matrix (have a reductive even part);

b) with a nondegenerate invariant odd supersymmetric bilinear form 
(have a reductive even part);

c) with a reductive even part (but not necessarily with a Cartan 
matrix or a nondegenerate form);

d) with nonreductive even part.

And the instances when the Weyl group appears 
are also numerous, here are a few examples (ordered historically):

A) as the set over which one performs the summation in the character 
formula (Bernstein--Leites character formula for atypical 
representations of $\fosp(2|2n)$ and $\fsl(1|n)$; or its conjectural 
generalization due to Penkov and Serganova);

B) as the group generated by ``reflections in simple roots'' 
(Skornyakov; Serganova; Egorov);

C) as the set of {\it neighboring} systems of simple roots.  The 
neighboring systems that are related by reflections in even roots can 
be, as well, considered equivalent under the action of the Weyl group 
of the even part of $\fg$; the other systems are related by 
reflections in odd roots whose product is not a priori defined  
(Leites and Shchepochkina call this set {\it skorpenser} in honor of 
Skornyakov, Penkov and Serganova who considered various definitions of 
{\it odd reflections}, see \cite{LSh});

D) as the group whose elements number the Schubert supercells (Manin).

How to define the reflection in the odd root $\alpha$ is not 
immediately clear.  For example, in one of the simplest cases of the 
general matrix algebra $\fgl(m|n)$ or its supertraceless subalgebra, 
for any odd root $\alpha$ we have $(\alpha, \alpha)=0$ whareas in the 
formula for the reflection one has to divide by this scalar product.

Though Skornyakov, Penkov and Serganova as well as Egorov and (in a 
different framework) Manin suggested several working definitions of 
the analog of the Weyl group, the very diversity of answers was one of 
the reasons why in this paper I try to reformulate the 
Chevalley theorem so as to avoid appealing to the notion of the Weyl 
group and make use of the root decomposition only.

Such formulation is formally applicable to any Lie algebra and Lie 
superalgebra provided one accordingly modifies the notion of ``the 
root decomposition''.  The corresponding statement is offered as a 
Conjecture.

To lighten presentation, already involving many cases, the 
description of invariant polynomials on $\fq$ and its relatives, 
$\fsq$ and $\fp\fsq$, will be given elsewhere together with the 
description of invariant polynomials on the Poisson superalgebra $\fpo$ 
and its relatives.

\ssec{0.3.  Earlier results} The first to completely describe the 
invariant polynomials was F.~Berezin who did it for a real form of 
$\fgl(m|n)$, see \cite{Be1}.  Simultaneously he started to consider the 
general case with V.~Kac but they split and published their results 
(obtained by distinct methods) separately \cite{Be2} and \cite{K2}.

Berezin's proof \cite{Be2} was based on analytical methods; its 
presentation is not user-friendly.  Hardly anybody really went through 
the proofs; statements, nevertheless, were a source of inspiration for 
several researchers.  Observe that Berezin only considered Lie 
superalgebras with root spaces of multiplicity one and with a 
nondegenerate invariant supersymmetric even bilinear form.

In the addition to the above, Kac \cite{K2} also described a rough 
structure of the algebra $I(\fg)$ of invariant polynomials on the Lie 
superalgebra $\fg$ of the same class that Berezin considered.  Kac's 
proof is rather lucid.

Berezin and Kac showed that there exists a polynomial $Q$ such that 
the localization of $I(\fg)$ with respect to the multiplicative system 
generated by $Q$ contains all $W(\fg_{\ev})$-invariant polynomials.

In \cite{S1}, \cite{S2} I described (a bit more explicitly than in 
\cite{Be1}, \cite{Be2} or \cite{K1}, \cite{K2}) invariant polynomials 
on Lie superalgebras of series $\fgl$, $\fsl$, $\fpe$, $\fspe$, $\fq$, 
$\fvect$, $\fsvect$, $\widetilde{\fsvect}$.

\ssec{0.4.  Main result and conjectures} Let $\fg$ be one of the 
following series or exceptional Lie superalgebras (in suggestive 
notations from the review \cite{L}):
$$
\fgl,\; \fsl,\; \fpsl,\; \fosp,\; \fpe,\; \fspe,\; 
\fosp_\alpha(4|2);\quad \fag_2,\; \fab_3;\quad \fvect,\; \fsvect,\; 
\widetilde{\fsvect}.\eqno{(0.1)}
$$
Let $\fh$ be the split Cartan subalgebra of $\fg$; let $R=R_+\cup R_-$ 
be the set of nonzero roots of $\fg$ divided into subsets of positive 
and negative ones.  Let
$$
\tilde R_+=\{\alpha \in R_+\mid -\alpha \in R_-\}.
$$

For each root $\alpha \in \tilde R_+$ denote by $\fg(\alpha)$ the Lie 
superlagebra generated by $\fh$ and the root superspaces $\fg_\alpha$ 
and $\fg_{-\alpha}$, i.e. $\fg(\alpha)=\fh\bigoplus\left 
(\mathop{\oplus}\limits_{k \in\Zee\setminus 0} \fg_{k\alpha}\right )$.  
This superalgebra plays an important role in various constructions.

In the cases considered in this paper it is isomorphic to the direct 
sum of $\fsl(2)$ or $\fosp(1|2)$ or $\fsl(1|1)$ and a commutative 
subalgebra of $\fh$. 

Let the image of $S(\fg(\alpha)^*)^{\fg(\alpha)}$ under the 
restriction homomorphism $S(\fg^*)\tto S(\fh^*)$ be denoted by 
$I^{\alpha}(\fh^*)$ and the image of $S(\fg^*)^{\fg}$ by $I(\fh^*)$.

\begin{Theorem}  {\em (Main Theorem.)} For the Lie superalgebras 
$(0.1)$ we have
$$
I(\fh^*)=\mathop{\cap}\limits_{\alpha\in \widetilde
R_+}I^{\alpha}(\fh^*).\eqno{(0.2)}
$$
If $\tilde R_+=\emptyset$, then $(0.2)$ should be read as 
$I(\fh^*)=S(\fh^*)^\fh$.  This statement is true for semisimple Lie 
algebras as well.
\end{Theorem}  

Observe that, except for $\fg=\fosp(1|2n)$, the algebra $I(\fh^*)$ is 
not noetherian whereas the localization of $I(\fh^*)$ with respect to 
the multiplicative system generated by $Q$ is noetherian.

Define $W(\fg)$ to be
$$
\text{the Weyl group of }\left\{\begin{array}{l} \fg_{\ev}\text{ if $\fg$ 
is of type $\fgl$, $\fsl$, $\fpsl$, $\fosp$, $\fpe$, $\fspe$, 
$\fosp_\alpha(4|2)$}\cr
\fg_{0}\text{ if $\fg$ 
is of type $\fvect$, $\fsvect$, $
\widetilde{\fsvect}$.}\end{array}\right.\eqno{(0.3)}
$$. 

\ssbegin{0.4.1}{Corollary} Let $Q=\mathop{\prod}\limits_{\alpha\in 
\widetilde R_{\od}}\alpha$.  Then any element of $S(\fh^*)^{W(\fg)}$ 
can be represented in the form $P/Q$, where $P\in I(\fh^*)$.
\end{Corollary}

\begin{Problem} Compute the Poincar\'e series of the algebras 
$I(\fh^*)$.  Observe that their description given in sec.  0.6 depends 
on the choice of the system of symple roots.  It is, perhaps, 
desirable, to give all such descriptions for all $W(\fg)$-inequivalent 
systems of simple roots, for the list of the latter, found by 
Serganova, see \cite{LSS}. 
\end{Problem}

\ssbegin{0.4.2}{Conjecture} Equality $(0.2)$ holds for all simple 
(and, moreover, related to them ``classical", see \cite{LSh}) Lie 
superlagebras.
\end{Conjecture}

\ssbegin{0.4.3}{Conjecture} Equality $(0.2)$ holds for all Lie algebra 
and Lie superlagebras provided we consider the {\em generalized weight 
decompositions} of Penkov--Serganova \cite{PS}.
\end{Conjecture}

Recall that if $\fg$ is a nilpotent Lie superalgebra and $V$ a 
$\fg$-module, then the {\em generalized weight decomposition} of $V$ 
with respect to $\fg$ is a presentation 
$V=\mathop{\oplus}\limits_{\mu\in\fg^{*}}V(\mu)$, where $\mu$ is such 
that $\mu([\fg_{\ev}, \fg_{\ev}])=0$ and $V(\mu)$ is the maximal 
$\fg$-submodule all the irreducible subquotients of which are 
isomorphic to $\Ind^{\fg}_{\fb}(\mu)$, where $\fb$ is the polarization 
for $\mu$, see \cite{K1}.

We suggest the reader to compare our approach with Shander's approach 
to invariant nonpolynomial {\it functions} on Lie superalgebras, cf.  
\cite{Sh}.

\ssec{0.5.  Chevalley's theorem for anti-invariant polynomials} It is 
well-known (\cite{Bu1}) that there exists a unique extention of the 
adjoint representation of a Lie algebra $\fg$ to representations in 
$S(\fg)$ and $U(\fg)$ and the canonical symmetrization $\omega: 
S(\fg)\tto U(\fg)$ is a $\fg$-module homomorphism.  For Lie 
superalgebras this supersymmetrization is given by the formula
$$
\omega(x_1, \dots , x_n)=\frac{1}{n}\mathop{\sum}\limits_{\sigma\in 
\fS_{n}}c(p(x), \sigma)x_{\sigma(1)} \dots 
x_{\sigma(n)},\eqno{(0.5.1)}
$$
where $p(x)$ is the vector of parities of $x=(x_{1}, \dots , x_{n})$ 
and $c(p(x), \sigma)$ is defined (see \cite{S1}) as follows.  
Let $a_{1}$, \dots , $a_{n}$ be the elements of the free 
supercommutative superalgebra.  Then
$$
c(p(x), \sigma)a_{1}
\dots a_{n}=a_{\sigma(1)} \dots a_{\sigma(n)}\text{ for any 
$\sigma\in\fS_{n}$}.\eqno{(0.5.2)}
$$

Consider also the continuation $\tilde \omega: S(\fg)\tto U(\fg)$ of 
the symmetrization $\omega: S(\fg_\ev)\tto U(\fg_\ev)$, where $S(\fg)$ 
is considered as $\Ind^\fg_{\fg_\ev}(S(\fg_\ev))$ and $U(\fg)$ is 
considered as $\fg$-module with the action
$$
x*u=(-1)^{p(x)}xu-(-1)^{p(x)p(u)}ux\text{ for $x\in\fg$, $u\in 
U(\fg)$}.\eqno{(0.5.3)}
$$

\begin{Theorem} Let $\fg$ be one of the Lie superalgebra from our list 
$(0.1)$, $\fh$ its split Cartan subalgebra, $W$ the Weyl group defined 
as in $(0.3)$ Then the restriciton homomorphism $S(\fg^*)\tto 
S(\fh^*)$ onto Cartan subalgebra induces an isomorphism
$$
S(\fg^*)^{\fg}\cong S(\fh^*)^{W}.
$$
\end{Theorem}  

\ssec{0.6. Summary: description of the algebra of invariant 
polynomials}  {}~{}

\par 
\underline{0.6.1.  $\fg=\fgl(n|m)$}. Let $\eps_1$, \dots , 
$\eps_n$, $\delta_1$, \dots , $\delta_m$ be the weights of the 
standard (identity) $\fgl(n|m)$-module in the standard basis.  The 
Weyl group is $W=\fS_n\times \fS_m$; it acts on the weights by 
separately permuting the $\eps_i$ and the $\delta_j$.  

We identify $S(\fh^*)$ with $\Cee[\eps_1, \dots , \eps_n; \delta_1, 
\dots , \delta_m]$.  Then
$$
I(\fh^*)=\{f\in \Cee[\eps_1, \dots , \eps_n,
\delta_1, \dots , \delta_m]^W\mid
\pderf{f}{\eps_i}+\pderf{f}{\delta_j}\in
(\eps_i-\delta_j)\}. 
$$
This is the algebra of supersymmetric polynomials.  For its generators 
we can take either power series
$$
\Delta_k=\sum \eps_i^k-\sum \delta_j^k
$$
or the coefficients of powers of $t$ in the power series expansion of the
rational function 
$$
F(t)=\frac{\prod(t-\delta_j)}{\prod(t-\eps_i)}.
$$
Denote $I(\fh^*)$ by $I_{n, m}$.

\underline{0.6.2. a) $\fg=\fsl(n|m)$}. $I(\fh^*)=I_{n, m}/(\eps_1+ \dots +
\eps_n-\delta_1-\dots -\delta_m)$.

\underline{0.6.2. b) $\fg=\fpsl(n|n)$}. Let 
$$
\tilde I_{n, n}=I_{n, n}\cap
\Cee[\eps_i-
\eps_j, \eps_i-\delta_j, \delta_i-\delta_j\text{ for 
$1\leq i,
j\leq n$}].
$$ 
Then $I(\fh^*)=\tilde I_{n, n}/(\eps_1+ \dots +
\eps_n-\delta_1-\dots -\delta_n)$.

\underline{0.6.3.  $\fg=\fosp(2m+1|2n)$, $m>0$}. We identify $S(\fh^*)$ 
with $\Cee[\eps_1, \dots , \eps_n, \delta_1, \dots , \delta_m]$ 
generated as algebra by the weights of the standard module.  The Weyl 
group is $W=(\fS_n\circ\Zee^n_2)\times (\fS_m\circ\Zee^m_2)$; it acts 
on the weights by separately permuting the $\eps_i$ and the $\delta_j$ 
and by changing the signs of the weights.
$$
I(\fh^*)=\{f\in \Cee[\eps_1, \dots , \eps_n,
\delta_1, \dots , \delta_m]^W\mid
\pderf{f}{\eps_i}\pm \pderf{f}{\delta_j}\in
(\eps_i\mp \delta_j)\}. 
$$
This is the algebra of supersymmetric polynomials $I_{n, m}$ in 
$\eps_1^2, \dots , \eps_n^2, \delta_1^2, \dots , 
\delta_m^2$.  For its generators we can take either power series
$$
\Delta_{2k}=\sum \eps_i^{2k}-\sum \delta_j^{2k}
$$
or the coefficients of powers of $t$ in the power series expansion of the
rational function 
$$
F(t)=\frac{\prod(t^2-\eps_i^2)}{\prod(t^2-\delta_j^2)}.
$$

\underline{0.6.4. $\fg=\fosp(1|2n)$}. This is a particular case: $F(t)$ is 
a polynomial.

\underline{0.6.5.  $\fg=\fosp(2m|2n)$, $m\geq 1$}.  We identify $S(\fh^*)$ 
with $\Cee[\eps_1, \dots , \eps_n, \delta_1, \dots , \delta_m]$, as in 0.6.3.

The Weyl group is $W=(\fS_m\circ\Zee^{m-1}_2)\times 
(\fS_n\circ\Zee^n_2)$; it acts on the weights by separately permuting 
the $\eps_i$ and the $\delta_j$ and $\Zee^{m-1}_2$ changes the 
signs of the weights $\eps_\mapsto \theta_i\eps_i$, 
where $\theta_i=\pm 1$ and $\prod \theta_i=1$.  In this case
$$
I(\fh^*)=\{f\in \Cee[\eps_1, \dots , \eps_m,
\delta_1, \dots , \delta_n]^W\mid
\pderf{f}{\eps_i}\pm \pderf{f}{\delta_j}\in
(\eps_i\mp \delta_j). 
$$
Any element of $I(\fh^*)$ can be expressed in the form
$$
f=f_0+\eps_1 \dots \eps_m\mathop{\prod}\limits_{i,
j}(\eps_i^2-\delta_j^2)\cdot f_1,
$$
where $f_0\in I_{n, m}(\eps_1^2, \dots , \eps_n^2, \delta_1^2, \dots , 
\delta_m^2)$ and can be expressed via the coefficients of the powers 
of $t$ in the power series expansion of the rational function
$$
F(t)=\frac{\prod(t^2-\delta_j^2)}{\prod(t^2-\eps_i^2)},
$$
and where $f_1\in\Cee[\eps_1, \dots , \eps_n, \delta_1, \dots
, \delta_m]^W$.  

\underline{0.6.6.  $\fg=\fosp_\alpha (4|2)$}. We identify $S(\fh^*)$ 
with the algebra $\Cee[\eps_1, \eps_2, \eps_3]$.

The Weyl group is
$W=\Zee_2\times \Zee_2\times \Zee_2$; it acts on the weights
by changes of their signs. Set
$$
\lambda_1=-(1+\alpha),\quad \lambda _2=1,\quad \lambda_3=\alpha.
$$ 
In this case for $\theta_i=\pm 1$ we have
$$
I(\fh^*)=\{f\in \Cee[\eps_1, \eps_2, 
\eps_3]^W\mid \theta_1\lambda_1\pderf{f}{\eps_1}+ 
\theta_2\lambda_2\pderf{f}{\eps_2}+ 
\theta_3\lambda_3\pderf{f}{\eps_3} \in (\theta_1\eps_1+ 
\theta_2\eps_2+\theta_3\eps_3)\}.
$$
Any element of $I(\fh^*)$ can be expressed in the form
$$
f=f_0+\prod(\eps_1\pm \eps_2\pm \eps_3) \cdot 
f_1,
$$
where $f_0\in \Cee[\frac{1}{\lambda_{1}}\eps_{1}^2+ 
\frac{1}{\lambda_{2}}\eps_{2}^2+ 
\frac{1}{\lambda_{3}}\eps_{3}^2]$ and 
$f_1\in\Cee[\eps_{1}^2, \eps_{2}^2, \eps_{3}^2]$.

\underline{0.6.7. $\fg=\fag_2$}. We identify $S(\fh^*)$ with the
algebra  $\Cee[\eps_1, \eps_2, \eps_3,
\delta]/(\eps_1+\eps_2+\eps_3)$. 

The Weyl group is $W=(\fS_3\circ\Zee_2)\times \Zee_2$; 
$\fS_3\circ\Zee_2$ acts on the $\eps_i$ by permutations and 
{\it simultaneous} changes their signs, the second factor, $\Zee_2$, 
changes the sign of the $\delta$.  In this case
$$
I(\fh^*)=\{f\in S(\fh^*)^W\mid
\pderf{f}{\eps_1}+\pderf{f}{\eps_2}+2\pderf{f}{\delta}
\in (\delta-\eps_1-\eps_2)\}.
$$
Any element of $I(\fh^*)$ can be expressed in the form
$$
f=f_0+\mathop{\prod}\limits_{1\leq i \leq 3}(\delta^2-\eps_i^2)\cdot 
f_1,
$$
where $f_0\in \Cee[3\delta^2-2(\eps_{1}^2+\eps_{2}^2+ 
\eps_{3}^2)]$ and $f_1\in S(\fh^*)^W$.

\underline{0.6.8.  $\fg=\fab_3$}. We identify $S(\fh^*)$ with 
$\Cee[\eps_1, \eps_2, \eps_3, \delta]$.

The Weyl group is $W=(\fS_3\circ\Zee_2^3)\times \Zee_2$; 
$\fS_3\circ\Zee_2^3$ acts on the $\eps_i$ by permutations and 
changes of their signs, the second factor, $\Zee_2$, changes the sign 
of the $\delta$.  In this case
$$
I(\fh^*)=\{f\in S(\fh^*)^W\mid 
\pderf{f}{\eps_1}+\pderf{f}{\eps_2}+ 
\pderf{f}{\eps_3}-3\pderf{f}{\delta} \in 
(\delta+\eps_1+\eps_2+\eps_3)\}.
$$
Any element of $I(\fh^*)$ can be expressed in the form
$$
f=f_0+\prod(\delta\pm \eps_1\pm \eps_2\pm 
\eps_3)\cdot f_1,
$$
where $f_0\in \Cee[L_2, L_6]$, $f_1\in S(\fh^*)^W$ and 
$$
\renewcommand{\arraystretch}{1.2}
\begin{array}{rl}
L_2&=3(\eps_{1}^2+\eps_{2}^2+\eps_{3}^2)-\delta^2\\
L_6&=\delta^6+\eps_{1}^6+\eps_{2}^6+\eps_{3}^6+ \\
&+(\eps_{1}-\eps_{2})^6+(\eps_{1}-\eps_{3})^6 
+(\eps_{2}-\eps_{3})^6+(\eps_{1}+\eps_{2})^6 +(\eps_{1}+\eps_{3})^6 
+(\eps_{2}+\eps_{3})^6-\\
&-\frac{1}{64}(\sum(\delta\pm \eps_1\pm \eps_2\pm
\eps_3)^6).
 \end{array}
$$
 
\underline{0.6.9.  $\fg=\fpe(n)$}. We identify $S(\fh^*)$ with the algebra 
$\Cee[\eps_1, \dots , \eps_n]$ generated by the weights 
of the standard $\fg$-module.

The Weyl group is $W=\fS_n$ acts on the $\eps_i$ by 
permutations.  In this case
$$
I(\fh^*)=\{f\in S(\fh^*)^W\mid
\pderf{f}{\eps_i}-\pderf{f}{\eps_j}
\in (\eps_i+\eps_j)\}.
$$
This is the algebra of projective Schur functions, for a system of its
generators we can take either the sums of powers
$$
\Delta_{2k+1}=\sum x_i^{2k+1}
$$
or the coefficients of the rational function
$$
F(t)=\frac{\prod(t+\eps_i)}{\prod(t-\eps_i)}.
$$  
Set $I_n=I(\fh^*)$

\underline{0.6.10. $\fg=\fspe(n)$}.  Then 
$I(\fh^*)=I_n/(\eps_1+\dots +\eps_n)$.

\underline{0.6.11.  $\fg=\fvect(0|n)$}. We identify $S(\fh^*)$ with the 
algebra $\Cee[\eps_1, \dots , \eps_n]$ generated by the weights of the 
$\fg_0$-module $\fg_{-1}$.

The Weyl group is $W=\fS_n$ acts on the $\eps_i$ by permutations.

In this case
$$
I(\fh^*)=\{f\in S(\fh^*)^W\mid \pderf{f}{\eps_i} \in (\eps_j)\text{ 
for any } i\neq j\}.
$$
Each element of $I(\fh^*)$ is of the form
$$
f=c+\eps_1 \dots \eps_ng,\text{ where } c\in\Cee, g\in
S(\fh^*)^W. 
$$

Set $J_n=I(\fh^*)$.

\underline{0.6.11.a) $\fg=\fsvect(0|n)$}.  Then $I(\fh^*)\cong 
J_n/(\eps_1+ \dots +\eps_n)$.

\underline{0.6.11.b) $\fg=\widetilde{\fsvect}(0|n)$}. Then the algebra
of invariant polinomials is the same as in the non-deformed case,
i.e., $I(\fh^*)\cong J_n/(\eps_1+ \dots +\eps_{2n})$.

\section*{\S 1. Preliminaries}

\ssbegin{1.1}{Proposition} {\em (\cite{S2})} Let $\fg$ be a finite 
dimensional Lie superalgebra, $V$ a finite dimensional $\fg$-module, 
$L\subset V$ a subsuperspace and $w_0\in W_{\ev}$ an element such that 
the map $\fg\times W\tto V$ given by the formula
$$
x, w \mapsto xw_0+w
$$
is surjective.  Then the restriction map $S(V^*)^\fg\tto S(W^*)$ is 
injective.
\end{Proposition}

\begin{proof} Fro notations see Appendix. Let $\rho : \fg_{\Lambda }\tto 
\fgl(V_{\Lambda})$ be a homomorphism.  Let $G_{\Lambda }$ be a 
connected and simply connected Lie group corresponding to 
$\fg_{\Lambda }$.  By the Lie theory (see \cite{OV}) there exists a 
unique homomorphism $\pi: G_{\Lambda }\tto GL(V_{\Lambda })$ such that 
the derivative of $\pi $ at the unit is equal to $\rho $.  Consider 
the manifold morphism
$$
\Phi : G_{\Lambda }\times L_{\Lambda }\tto V_{\Lambda }, \quad \Phi 
(g, w)=gw.
$$
It is easy to calculate the derivative of this map.  It is equal to
$$
(D\Phi)(1, w_{0})(\tilde{x}, \tilde{w})=\rho 
(\tilde{x})w_{0}+\tilde{w}.\eqno{(*)}
$$
 
Let $f\in S(V^{*})^{\fg}$ and $i(f)=0$.  Then $\theta 
(f)|_{L_{\Lambda}}=0$.  By Lemma A6 $\theta (f)\in \Lambda \otimes 
S(V^{*}_{\Lambda})$ is $\fg_{\Lambda}$-invariant.  By Corollary 3 of 
Proposition3.6.13 from \cite{Bu2} $\theta (f)$ is invariant with 
respect to the natural $G_{\Lambda }$-action on $\Lambda \otimes 
S(V^{*}_{\Lambda})$.
 
 Consider $\theta (f)$ as the usual polynomial mapping 
 $V_{\Lambda}\tto \Lambda$.  It follows from $(*)$ that $\Phi $ is 
 submersive at $(1, w_{0})$.  Therefore, the image of $\Phi $ contains 
 an open neighborhood of $w_{0}$ in $V_{\Lambda }$.  Hence, 
 $G_{\Lambda }L_{\Lambda }$ contains an open neighborhood of $w_{0}$; 
 hence, $\theta (f)=0$.  Lemma A6 shows that $f=0$.  \end{proof}

\begin{Corollary} If $\fg$ is a finite dimensional Lie superalgebra, 
$\fh$ its its Cartan subalgebra, then the restriction map 
$S(\fg^*)^\fg\tto S(\fh^*)$ is injective.  
\end{Corollary}

\begin{proof} It suffices to demonstrate that  the map
 $\fg\times \fh\tto \fg$ given by the formula
$$
x, h \mapsto xh_0+h
$$
is surjective.

Let $\fg=\mathop{\oplus}\limits_{\alpha\in R} \fg(\alpha)$ be the 
generalized weight decomposition of $\fg$ relative $\fh$, see sec. 0.4.3.

Select $h_0\in\fh_\ev$ such that $\alpha(h_0)\neq 0$ for all 
$\alpha\in R$.  Then the kernel of the above map consists of the pairs 
$(x, h)$ such that $xh_0+h=0$.  Since $\alpha(h_0)\neq 0$ for all 
$\alpha\in R$, it folloows that $x\in \fh$ and the dimension of the 
kernel is equal to that of $\fh$.  Hence, the map $\fg\times \fh\tto 
\fg$ is surjective.
\end{proof} 
	
\ssbegin{1.2}{Proposition} Let $\fg$ be a finite dimensional Lie 
superalgebra such that the Lie algebra $\fg_{\ev}=\fh$ is commutative; 
$\fg_\od=\Span(u, v)$ and the following relations hold:
$$ 
\renewcommand{\arraystretch}{1.2}
\begin{array}{l}
{}[h, v]=-\alpha(h)v,\; [h, u]=\alpha(h)u \text{ for any $h\in \fh$
and $\alpha\in\fh^*$};\\
{}[u, v]=h_{\alpha} \text{ and $\alpha(h_{\alpha})=0$.}\\
\end{array}
$$
Let 
$$
I^\alpha(\fh^*)=\{f\in S(\fh^*)\mid D_{h_{\alpha}}f\in(\alpha)\},
$$
where $D_{h_{\alpha}}f$ is the derivative of $f$ in the diriction of 
$h_\alpha$.

Then the restriciton homomorphism $S(\fg^*)\tto S(\fh^*)$ induces an 
isomorphism $S(\fg^*)^\fg\cong I^\alpha(\fh^*)$.
\end{Proposition}

\begin{proof} Since $\fh$ is the Cartan subalgebra, the restriction 
homomorphism is injective by Proposition 1.1.  Let us prove that its 
image coincides with $I^\alpha(\fh^*)$.

Let $F\in S(\fg^*)^\fg$. By $\fh$-invariance $F$ should be of the form 
$$
F=f+gv^*u^*, \text{ where $v^*$ and $u^*$ is the left dual basis to 
$v$ and $u$}.
$$
It is easy to verify that in $S(\fg^*)$ the following relations hold:
$$
\renewcommand{\arraystretch}{1.2}
\begin{array}{l}
v\cdot l=-l(h_{\alpha})\cdot u^*,\; 
u\cdot l=-l(h_{\alpha})\cdot v^*,\; \text{ for $l\in\fh^*$}\\
v\cdot v^*=\alpha;\; v\cdot u^*=0;\; u\cdot u^*=-\alpha;\; u\cdot 
v^*=0.
\end{array}
$$
Therefore, $v\cdot F=-D_{h_{\alpha}}f\cdot u^*+\alpha gu^*=0$.  Hence, 
$D_{h_{\alpha}}f=\alpha g$ and $f\in I^\alpha(\fh^*)$.

Convesely, if $f\in I^\alpha(\fh^*)$, then it is easy to verify that
$$
F=f+\frac{1}{\alpha}D_{h_{\alpha}}fv^*u^*\in S(\fg^*)^\fg.
$$
\end{proof}

\ssbegin{1.3}{Proposition} Let $\fg$ be a finite dimensional Lie 
superalgebra such that the Lie algebra $\fg_{\ev}=\fh$ is commutative; 
$\fg_\od=\Span(v_1, u_1, v_2, u_2)$ and the following relations hold:
$$ 
\renewcommand{\arraystretch}{1.2}
\begin{array}{l}
{}[h, v_i]=-\alpha(h)v_i,\; 
[h, u_i]=\alpha(h)u_i \text{ for any $h\in \fh$, $i=1,
2$ and $\alpha\in\fh^*$};\\ 
{}[u_i, u_j]=[v_i, v_j]=0;\\
{}[u_i, v_j]=0 \text{ for $i\neq j;\; [u_i, v_i]= h_i$ (it is possible that
$h_1=h_2$);}\\
\alpha([\fg_\od, \fg_\od])=0.
\end{array}
$$
Let 
$$
I^\alpha(\fh^*)=\{f\in S(\fh^*)\mid D_{h_{1}}f\in (\alpha),\;
D_{h_{2}}f\in (\alpha), \; D_{h_{1}}D_{h_{2}}f\in (\alpha^2)\}, 
$$
where $D_{h}f$ is the derivative of $f$ in the direction of 
$h$.

Then the restriciton homomorphism $S(\fg^*)\tto S(\fh^*)$ induces an 
isomorphism $S(\fg^*)^\fg\cong I^\alpha(\fh^*)$.
\end{Proposition}

\begin{proof} As in the proof of Proposition 1.2, let $F\in 
S(\fg^*)^\fg$.  By $\fh$-invariance, $F$ should be of the form
$$
F=f+gv_1^*u_1^*+sv_2^*u_2^*+rv_1^*u_1^*v_2^*u_2^*, 
$$
where the $v_i^*$ and $u_i^*$ constitute the left dual basis to the 
$v_i$, $u_i$.

The condition $v_{1}F=0$ implies that
$$
D_{h_{1}}f-\alpha g=0,\quad D_{h_{1}}s-\alpha r=0,  
$$ 
and similarly
$$
D_{h_{2}}f-\alpha s=0,\; D_{h_{2}}s-\alpha r=0.  
$$ 
These conditons imply that
$$
D_{h_{1}}f\in (\alpha),\quad 
D_{h_{2}}f\in (\alpha), \quad D_{h_{1}}D_{h_{2}}f\in (\alpha^2).
$$

Convesely, if these conditions hold, then by setting
$$
F=f+\frac{1}{\alpha}D_{h_{1}}fv_1^*u_1^*+
\frac{1}{\alpha}D_{h_{2}}fv_2^*u_2^*+ 
\frac{1}{\alpha^2}D_{h_{1}}D_{h_{2}}fv_1^*u_1^*v_2^*u_2^*
$$
we get an element from $S(\fg^*)^\fg$ whose restriction on $\fh^*$ is 
equal to $f$.  \end{proof}

\begin{Remark} Under conditions of Proposition 2.2, if $h_1=h_2$, then 
not every element from $S(\fg^*)^\fg$ can be obtained as a linear 
combination of invariant polynomials of the form $\str(\rho(x)^k)$, 
where $\rho$ is any finite dimensional representation of $\fg$.

Indeed, let $h_1=h_2=h$.  Then there are two types of representations 
of $\fg$:

1) the 1-dimensional ones, determined by linear forms 
$\lambda\in\fh^*$ such that $\lambda(h)=0$ and

2) $(4, 4)$-dimensional ones for which $\lambda(h)\neq 0$; these
representations are of the form 
$$
T_\lambda=\ind^\fg_{\fh\oplus\Span(u_1, u_2)}(\lambda).
$$

The character of a representation of the first type is equal (up to
$\eps$) to $e^{\lambda}$; these characters generate a subalgebra of
invariants isomorphic to $\Cee[\fh^\perp]$, where
$$
\fh^\perp=\{\lambda\in\fh^*\mid \lambda(h)=0\}.
$$

Further, it is not difficult to verify that each homogeneous component 
of the supercharacter in the sence of Kac considered as a formal power 
series is an invariant polynomial.  (Recall that $\sch V$ is a linear 
functional on $U(\fg)$ defined as $\sch V(u)=\str_{V}(u)$; here we 
identify $U(\fg)^{*}$ with the algebra of formal power series.)  
Actually, this component {\it is} $\str(\rho(x)^k)$, where $\rho$ is 
the corresponding representation of $\fg$.

Therefore, the homogeneous components of supercharacters of 
$T_\lambda$ are invariant polynomials.  But the supercharacter of 
$T_\lambda$ is equal to $(1-e^{-\alpha})^2e^{\lambda}$ and the linear 
span of such supercharacters is equal to $\alpha^2\cdot \varphi$, 
where $\varphi\in\Cee[\fh^*]$.  Thus, the general form of the 
restriction of an invariant polynomial equal to the linear combination 
of supertraces is of the form
$$
\psi+ \alpha^2\cdot \varphi,\text{ where $\varphi\in\Cee[\fh^*]$, 
$\psi\in\Cee[\fh^\perp]$}.
$$

On the other hand, by Proposition 1.3 the element $\alpha\cdot l$, 
where $l\in\fh^*$ and $l(h)\neq 0$, is a restriction of an invariant 
and does not have the above described form.  Therefore, $\alpha\cdot 
l$ does not belong to the subalgebra generated by supertraces of 
finite dimensional representations.
\end{Remark}

\ssbegin{1.4}{Proposition} Let $\fg=\mathop{\oplus}\limits_{-1\leq i\leq
N}\fg_i$ be a consistent grading of a $\Zee$-graded Lie superalgebra;
moreover, let $\fg_0$ be reductive. Let $W=W(\fg_0)$ be the  Weyl group of
$\fg_0$ and $\alpha_1$, \dots , $\alpha_n$ be the weights of the
$\fg_0$-module $\fg_{-1}$.   

Then any polynomial of the form $\alpha_1\dots\alpha_nf$ on the
Cartan subalgebra $\fh$ of $\fg_0$, where $f\in \Cee[\fh^*]^W$,
is a restriction of an invariant.
\end{Proposition}

\begin{proof} Let $\lambda$ be a highest weight for $\fg_0$, 
$L^{\lambda}$ the corresponding irreducible $\fg_0$-module (with even 
heighest weight vector).  Let $\fg_iL^{\lambda}=0$ for $i>0$; denote 
$\fg_{\geq 0}=\mathop{\oplus}\limits_{i\geq 1}\fg_i$ and set 
$V^{\lambda}=\ind^\fg_{\fg_{\geq 0}}(L^{\lambda})$.  Then
$$
\sch V^{\lambda}=\mathop{\prod}\limits_{1\leq i\leq n} 
(1-e^{\alpha_{i}})\ch L^{\lambda}.
$$
According to \cite{B} there exists a linear combination $\sum 
c_\lambda \ch L^{\lambda}$ of functions $\ch L^{\lambda}$ such that 
the homogeneous component of the least degree of this combination is 
equal to $f$.  But then the homogeneous component of the least degree 
of $\mathop{\prod}\limits_{1\leq i\leq n}(1-e^{\alpha_{i}})\sum 
c_\lambda \ch L^{\lambda}$ is equal to $\alpha_1\dots\alpha_nf$.
\end{proof}

\ssbegin{1.5}{Proposition} Let $\fg=\fg_\ev\oplus\fg_\od$ be a Lie 
superalgebra with reductive $\fg_\ev$.  Let $W=W(\fg_\ev)$ be the Weyl 
group of $\fg_\ev$ and $\alpha_1$, \dots , $\alpha_n$ be the weights 
of the $\fg_\ev$-module $\fg_\od$.

Then any polynomial of the form $\alpha_1\dots\alpha_nf$ on the
Cartan subalgebra $\fh$ of $\fg_\ev$, where $f\in \Cee[\fh^*]^W$,
is a restriction of an invariant.
\end{Proposition}

Proof is similar to that of Proposition 1.4.

\ssbegin{1.6}{Proposition} Let $\fg=\fsl(2)$, $L(n)$ the finite 
dimensional irreducible $\fg$-module with highest weight $n\varphi$, 
where $\varphi$ is the fundamental weight.

Then for any $k\in\Nee$ and any $n_1, \dots , n_{k+1}\in\Nee$ there 
exist $c_1, \dots , c_{k+1}\in\Cee$ such that the homogeneous 
component of the least degree of the linear combination of characters 
$\sum c_i \ch L^{n_{i}}$ under the restriction on the Cartan subalgebra 
becomes equal to $\varphi^{2k}$.
\end{Proposition}

\begin{proof}Since
$$
\ch L^{n_{i}}=\frac{\sinh(n_i+1)\varphi}{\sinh \varphi},
$$
it follows that it suffices to prove the statement for a linear 
combination $\sum c_i \sinh (n_i+1)\varphi$.  By equating the 
coefficients of the powers 1, 3, \dots , $2k-1$ to zero and equating 
the coefficient of the power $2k+1$ to 1 we get the system
$$
\begin{array}{rl}
c_1(n_1+1)+ \dots + c_{k+1}(n_{k+1}+1)&=0\\
c_1(n_1+1)^3+ \dots + c_{k+1}(n_{k+1}+1)^3&=0\\
\dotfill&\dotfill\\
c_1(n_1+1)^{2k-1}+ \dots + c_{k+1}(n_{k+1}+1)^{2k-1}&=0\\
c_1(n_1+1)^{2k+1}\dots + c_{k+1}(n_{k+1}+1)^{2k+1}&=(2k+1)!\\
\end{array}
$$
If the $n_i$ are pairwise distinct, the determinant of the system is 
nonzero.
\end{proof}

\begin{Remark} Proposition 1.6 implies that
$$
\mathop{\sum}\limits_{1\leq i\leq k+1} c_i \sinh (n_ix)\varphi= 
x^{2k+1}+\dots
$$
By differentiating this identity with respect to $x$ we see that there 
exist $c_1'$, \dots $c_{k+1}'$ such that
$$
\mathop{\sum}\limits_{1\leq i\leq k+1} c_i' \cosh 
(n_ix)\varphi=x^{2k}+\dots .
$$
\end{Remark}

\ssec{1.7.  The algebra of dual numbers and modules over associative 
superalgebras} Let $A$ be an associative superalgebra.  Consider the 
algebra (with any parity ignored)
$$
A[\eps]=\{a+b\eps\mid a, b\in A, \; \eps^2=1, \; \eps
a=(-1)^{p(a)}a\eps\}. 
$$
Observe that our algebra $A[\eps]$ is almost the same as the {\it 
super}algebra $Q(A)$ considered in \cite{BL}, but (1) we ignore 
parity, (2) in $Q(A)$ a different relation $\eps^2=-1$, should be 
satisfied.

\ssbegin{1.7.1}{Lemma} The category $\Mods$ of $A[\eps]$-modules is 
isomorphic to the category $\SMods$ of graded $A$-modules the 
morphisms of the latter being the purely even homomorphisms.
\end{Lemma}

\begin{proof} Let us construct the functors
$$
F:\Mods\tto \SMods\text{ and } G:\SMods \tto \Mods.
$$
Let $V\in\Ob \Mods$, then $\eps\in\End(V)$ and since 
$\eps^2=1$, it follows that $F(V)=V_\ev\oplus V_\od$, where
$$
V_\ev=\{v\in V\mid \eps v=v\}\text{ and } 
V_\od=\{v\in V\mid \eps v=-v\}.
$$
Thus, $V$ may be considered as an object from $\SMods$.  It is an easy 
and routine chech with the help of the commutation relations for 
$\eps$ that if $f: V\tto W$ is an $A[\eps]$-module homomorphism, then 
$f$ may be considered as an even homomorphism (morphism) of 
$A$-modules.

Every superspace $V=V_\ev\oplus V_\od\in\Ob \SMods$ has a fixed map 
(not a morphism), namely the parity map $J$ such that 
$J(v)=(-1)^{p(v)}v$, see \cite{Be2}.  Sending $\eps$ to $J$ we get an 
object from $\Mods$.

To verify that $FG$ and $GF$ are identity functors is trivial.
\end{proof}

\ssbegin{1.7.2}{Lemma} Let $A$ be an associative superalgebra and 
$L(A)$ the corresponding Lie algebra.  Set
$$
a*b=(-1)^{p(a)}ab-(-1)^{p(a)p(b)}ba.
$$
With respect to the action $*$ the space $A$ is an $L(A)$-module.  
Moreover, let $V$ be a graded $A$-module.  The function $a\mapsto 
\tr_V(a)$ is an $L(A)$-invariant function on $A$.  
\end{Lemma}

\begin{proof} Consider $A[\eps]$ as a superalgebra with $p(\eps)=\ev$.  
(Observe that in \cite{BL} $p(\eps)=\od$.)  Let us calculate $[a, \eps 
b]$, where $a, b\in A$ and $[\cdot , \cdot]$ denotes the 
supercommutator.  We have
$$
[a, \eps b]=a\cdot \eps b-(-1)^{p(a)p(b)}\eps ba=
\eps ((-1)^{p(a)}ab-(-1)^{p(a)p(b)}ba),
$$
which leads to the above formula for $a*b$.

Further on, $V$ is a graded $A[\eps]$-module (since $J$ is an 
even map); hence by the characteristic property of the supertrace 
$\str [a, \eps b]=0$; hence,
$$
\tr(a*b)=\str(\eps a*b)=\str[a, \eps b]=0
$$
(because if $a\in A$, then $\str(a)=\tr(\eps a)$ and
$\tr(a)=\str(\eps a)$).
\end{proof}

\begin{Remark} We are forced to introduce $A[\eps]$ in its 
latter interpretation because we wish to consider not only the 
intertwining operators that become scalars on the irreducible modules 
but may become a scalar multiple of $J$.  This train of thought leads 
us to the notion of anticenter considered elsewhere.
\end{Remark} 

\section*{\S 2. Invariant polynomials on Lie superlagebras}

Let $\fg$ be one of the Lie superalgebras 
$$
\fgl,\; \fsl,\; \fpsl,\; \fosp,\; \fpe,\; \fspe,\; \fosp_\alpha(4|2);\quad
\fag_2,\; \fab_3;\quad \fvect,\; \fsvect,\; \widetilde{\fsvect}.\eqno{(0.1)}
$$
Let $\fh$ be the split Cartan subalgebra of $\fg$; let $R=R^+\cup R^{-}$ 
be the set of nonzero roots of $\fg$ divided into subsets of positive 
and negative ones.  Let $R_{0}$ and $R_{1}$ the subsets of 
even and odd roots; each set can be divided into the subsets 
$R_{0}^\pm$, $R_{1}^\pm$ of positive and negative roots.  We further 
set:
$$
\tilde R^+=\{\alpha \in R^+\mid -\alpha \in R^-\}\quad 
\tilde R_{1}^+=\{\alpha \in R_{1}^+\mid -\alpha \in R_{1}^-;\;2\alpha 
\not\in R_{1}^+\}.\eqno{(2.1)}
$$
Let $\fg=\fh\oplus_{\alpha\in R}\fg^\alpha$ be the weight 
decomposition; for $\alpha\in \tilde R_{1}^+$ set 
$$
\nu(\alpha)=\min(\dim\fg^\alpha, \dim \fg^{-\alpha}).\eqno{(2.2)}
$$
For the Lie superalgebras considered $\nu(\alpha)=1$ or 2. We set
$$
I^{\alpha}(\fh^{*})=\left\{\renewcommand{\arraystretch}{1.2}
\begin{array}{l}
\{f\in S(\fh^{*})\mid D_{h}f\in(\alpha)\text{ 
for }h\in [\fg_{\od}^{-\alpha}, \fg_{\od}^{\alpha}]\}\; \text{ if 
$\nu(\alpha)=1$} \\
\{f\in S(\fh^{*})\mid D_{h}f\in(\alpha)\text{ 
for }h\in [\fg_{\od}^{-\alpha}, \fg_{\od}^{\alpha}]\\
\text{and $D_{h_{1}}D_{h_{2}}f\in (\alpha^2)$ 
for $h_{1}$, $h_{2}$ that generate } [\fg_{\od}^{-\alpha}, 
\fg_{\od}^{\alpha}]\}\; \text{ if $\nu(\alpha)=2$} .
\end{array}\right .
$$
Let $W$ be the Weyl group of $\fg_{\ev}$ if $\fg_{\ev}$ is reductive 
and the Weyl group of $\fg_{0}$ if $\fg$ is a vectorial Lie 
superalgebra (with nonreductive $\fg_{\ev}$) in the standard 
$\Zee$-grading.

\begin{Theorem} The restriction homomorphism $S(\fg^*)\tto 
S(\fh^*)$ induces an isomorphism of $S(\fg^*)^\fg$ with the algebra
$$
I(\fh^{*})=\{f\in S(\fh^{*})^{W}\mid f\in I^{\alpha}(\fh^{*})\text{ 
for any }\alpha \in \tilde R_{1}^+\}.
$$
\end{Theorem}

\begin{proof} Let us prove that the image $f$ of any invariant 
polynomial $F$ belongs to $I(\fh^{*})$. First, observe that $f\in 
S(\fh^{*})^{W}$.

Define: $\fg_{\alpha}$ to be the Lie subsuperalgebra of $\fg$ 
generated by the root vectors of weights proportional to $\alpha$.  
Clearly, if $\fg_\ev$ is reductive, $\fg_{\alpha}$ is isomorphic to 
one of the Lie (super)algebras from Proposition 1.2 if 
$\nu(\alpha)=1$ or  one of the Lie (super)algebras from Proposition 
1.3 if $\nu(\alpha)=2$ and $\fg=\fsl(2|2)$, $\fpsl(2|2)$ or $\fspe(4)$.

Observe that the restriction of 
$F$ onto $\fg_{\alpha}$ belongs to $S(\fg_{\alpha}^{*})^{\fg_{\alpha}}$.
From Propositions 1.2 and 1.3 we deduce that $f\in I^{\alpha}(\fh^{*})$.

If $\fg$ is of vectorial type, then $\nu(\alpha)=1$ and by selecting 
$u\in\fg^{-\alpha}$, $v\in\fg^{\alpha}$ we obtain a subalgebra 
$\fb=\fh\oplus\Span(u, v)$ satisfying the conditions of Proposition 1.2. 
Hence, in this case $f\in I^{\alpha}(\fh^{*})$ as well. 

Moreover, Proposition 1.1 shows that the restriction homomorphism is 
injective. Let us prove now that each element from $I(\fh^{*})$ can be 
extended to an element from $S(\fg^{*})^{\fg}$. The proof will be 
caried out separately for each of the above Lie superalgebras.

\underline{$\fg=\fgl(n|m)$}. (cf. \cite{B}) On the Cartan subalgebra, 
consider the formal power series $\Cee[[\fh^*]]$. For $h\in\fh$, 
$\lambda\in\fh^*$ and $e^{\lambda}\in\Cee[[\fh^*]]$ set
$$
D_{h}e^{\lambda}=\lambda(h)e^{\lambda}.
$$
Define
$$
\renewcommand{\arraystretch}{1.2}
\begin{array}{l}
J(\fh)^*=\{f\in\Cee[[\fh^*]]^{W}\mid f\text{ is a linear combination 
of the $e^{\lambda}$, }\\
\text{ where $\lambda$ is a weight of a 
representation of $\fg_{\ev}$ and}\\
D_{h}f\in(\alpha)\text{ for any $\alpha\in R_{\od}$ and $h\in 
[\fg_{\od}^{-\alpha}, \fg_{\od}^{\alpha}]$}.
\end{array}
$$

Let us prove that any element from $J(\fh^{*})$ is a linear 
combination of supercharacters of finite dimensional representations.  
Indeed, since the homogeneous components of supercharacters are 
invariant polynomials, it follows that any linear combination of them 
belongs to $J(\fh)^*$.  Let $(\fg_{\ev})_{s}$ be the semisimple part 
of $\fg_{\ev}$ and $\langle \cdot , \cdot\rangle$ the Killing form on 
$(\fg_{\ev})_{s}$.

For $\lambda\in\fh^{*}$ let $\lambda_{s}$ be its restriction onto 
$\fh\cap (\fg_{\ev})_{s}$; we set
$$
|\lambda|^2=\langle \lambda_{s}, \lambda_{s}\rangle.
$$
For $f=\sum c_{\nu}e^\nu$ set $r_{f}=\max_{\nu\neq 0}|\nu_{s}|$ and 
induct on $r_{f}$.

If $r_{f}=0$, then $\nu_{s}=0$ for any $\nu$.  Let $\alpha\in R_{\od}$ 
and $h\in [\fg_{\od}^{-\alpha}, \fg_{\od}^{\alpha}]$.  If 
$\nu(h)=0$, then $\nu_{s}=0$ implies that $\nu$ is 
proportional to the supertrace, the supercharacter of a 
one-dimensional representarion. Therefore, we may assume that 
$\nu(h)\neq 0$ for any $\nu$. 

Consider the restriction of $D_{h}f$ onto $[\fg_{\od}^{-\alpha}, 
\fg_{\od}^{\alpha}]$.  Since $\alpha(h)=0$, it follows that 
$D_{h}f=\sum c_{\nu}\nu(h)e^\nu=0$.  Since the exponents are linearly 
independent, $c_{\nu}\nu(h)=0$ for all $\nu$.  And since $\nu(h)\neq 
0$, it follows that each $c_{\nu}=0$; hence, $f=0$.  Thus, if 
$|\nu_{s}|=0$, then $f=\sum c_{\nu} e^\nu$, where each $\nu$ is 
proportional to the supertrace.  Let $r_{f}>0$ and $f=\sum 
c_{\nu}e^\nu$.

Consider the difference $f-\sum c_{\nu}\sch L^\nu$, where $L^\nu$ is 
therreducible highest weight module with the highest weight $\nu$ and 
the sum runs over $\nu$ such that $|\nu_{s}|=r_{f}$. Since the 
difference is $W$-invariant, it contains summands $e^\mu$ with 
$|\mu_{s}|<r_{f}$ by a property of representations of semisimple Lie 
algebras (see \cite{B}) and the inductive hypothesis applies. 

Now, let $P\in I(\fh^{*})$; consider $P$ as a polynomial in weights 
of the standard module, i.e., $P=P(\eps_{1}, \dots, \eps_{n}; 
\delta_{1}, \dots , \delta_{m})$ and let
$$
f=P(e^{\eps_{1}}-1, \dots, e^{\eps_{n}}-1; 
e^{\delta_{1}}-1, \dots , e^{\delta_{m}}-1).
$$
Let us verify that $f\in J(\fh^{*})$. Clearly, $f$ is $W$-invariant. 
Let $\alpha=\eps_{i}-\delta_{j}$ and $h\in [\fg_{\od}^{-\alpha}, 
\fg_{\od}^{\alpha}]$; then the condition $D_{h}f\in (\alpha)$ is 
equivalent to the fact that
$$
\renewcommand{\arraystretch}{1.4}
\begin{array}{l}
\pderf{f}{\eps_{i}}+\pderf{f}{\delta_{j}}\in(\eps_{i}-\delta_{j})
\Longleftrightarrow 
\pderf{P}{\eps_{i}}e^{\eps_{i}}+\pderf{P}{\delta_{j}}e^{\delta_{j}}\in ( 
e^{\eps_{i}}-e^{\delta_{j}})\Longleftrightarrow \\
e^{\delta_{j}}\left( \pderf{P}{\eps_{i}}e^{\eps_{i}-\delta_{j}}+
\pderf{P}{\delta_{j}}\right) \in ( 
e^{\eps_{i}-\delta_{j}}-1)e^{\delta_{j}}\Longleftrightarrow 
\left( \pderf{P}{\eps_{i}}e^{\eps_{i}-\delta_{j}}+
\pderf{P}{\delta_{j}}\right )\in (e^{\eps_{i}-\delta_{j}}-1).
\end{array}
$$
Since 
$\pderf{P}{\eps_{i}}+\pderf{P}{\delta_{j}}\in(\eps_{i}-\delta_{j})$, 
the last membership implies that 
$\pderf{f}{\eps_{i}}+\pderf{f}{\delta_{j}}\in(\eps_{i}-\delta_{j})$.  
Since $f\in J(\fh^{*})$, any homogeneous component of $f$ is a 
restriction of an invariant; but $P$ is an homogeneous component of 
$f$ of the least degree, hence, it is a restriction of an invariant.
\end{proof}

\underline{$\fg=\fsl(n|m)$, $n\neq m$}.  Let $\fh$ the Cartan 
subalgebra of $\fgl(n|m)$; let $\tilde\fh= \fh\cap \fg$.  Let $f\in 
I(\tilde\fh^*)$; set $g(h)=f\left(h-\frac{\str h}{n-m}1_{n+m}\right)$.  
Then
$$
(D_{h_{\alpha}}g)(h)=(D_{h_{\alpha}}f)(h-\frac{\str 
h}{n-m}1_{n+m})\in \alpha(h-\frac{\str 
h}{n-m}1_{n+m})=\alpha(h).
$$
Hence, $g\in I(\tilde\fh^*)$ and by the above is a restriction of an 
invariant.  Hence, $f$ is also a restriction of an invariant.

\underline{$\fg=\fsl(n|n)$, $n\neq 2$}. Let us show that the 
restriction map $I(\fh^{*})\tto I(\tilde\fh^{*})$ is surjective. Let 
us identify $S(\fh^{*})$ with $\Cee[\eps_{1}, \dots, \eps_{n}; 
\delta_{1}, \dots , \delta_{m}]$, where $\eps_{1}, \dots, \eps_{n}; 
\delta_{1}, \dots , \delta_{m}$ are the weights of the identity  
representation. It is not difficult to show (see\cite{P}) that 
$$
I(\fh^{*})\cong\Cee[s_{1}, s_{2}, \dots  ],\text{ where 
$s_{k}=\sum\eps_{i}^k-\sum\delta_{j}^k$ are supertraces of $\id$.}
$$
Let $\tilde\eps_{i}$ and $\tilde\delta_{j}$ be the images of 
$\eps_{i}$ and $\delta_{j}$ in $S(\tilde\fh^{*}))$. Since the 
$\tilde\eps_{i}$ and $\tilde\delta_{j}$ for $j<n$ are algebraically 
independent, we deduce, as above, that 
$$
I(\tilde\fh^{*})\subset\Cee[\sigma_{1}, \sigma_{2}, \dots  ],\text{ where 
$\sigma_{k}=\sum\tilde\eps_{i}^k-\sum\tilde\delta_{j}^k$.}
$$
If $f\in I(\tilde\fh^{*})$, then $f=F(\sigma_{1}, \sigma_{2}, \dots 
)$.  Consider $d=f-F(\sigma_{1}-\sigma_{1}, \sigma_{2}-\sigma_{1}^2, 
\dots )$.  Clearly $d=0$ at $\sigma_{1}=0$; hence, $d$ is divisible by 
$\sigma_{1}$. Since $F(\sigma_{1}-\sigma_{1}, \sigma_{2}-\sigma_{1}^2, 
\dots )$ is the image of an element from $I(\fh^{*})$, we may assume 
that $f$ is divisible by $\sigma_{1}$, i.e., $f=\sigma_{1}g$. 

For $\alpha=\eps_{n}-\delta_{n}$ the condition 
$D_{h_{\alpha}}f\in(\alpha)$ is equivalent to the fact that the 
restriction of $f$ onto $\ker \alpha$ is invariant with respect to 
translations by the vectors $h_{\alpha}$. Since 
$\sigma_{1}(h_{\alpha})\neq 0$, any element from $\ker \alpha$ can be 
expressed in the form $h+th_{\alpha}$, where $h\in\ker \alpha\cap 
\ker\sigma_{1}$. Therefore,
$$
f(h+th_{\alpha})=\sigma_{1}(h+th_{\alpha})g(h+th_{\alpha})=
t\sigma_{1}(h_{\alpha})g(h+th_{\alpha})=f(h)=0.
$$
Therefore, $g(h+th_{\alpha})=$ and the restriction of $f$ onto $\ker 
\alpha$ is equal to 0.  Hence, $f$ is divisible by $\alpha$; hence, 
from the $W$-symmetry and the fact that the linear functions $w\alpha$ 
are pair-wise coprime for $n>2$ we deduce that
$$
f=\prod_{\alpha\in R_{\od}^+}\alpha\cdot \tilde\varphi,\text{ where 
}\tilde\varphi,\in\Cee[\tilde\fh^{*}]^W.
$$
It is clear that $f$ is a restriction of an element of the form
$$
\prod_{\alpha\in R_{\od}^+}\alpha\cdot \varphi,\text{ where 
}\varphi,\in\Cee[\fh^{*}]^W.
$$
Proposition 1.4 applied to 
$\fg=\fgl(n|n)=\fg_{-}\oplus\fh\oplus\fg_{+}$, where $\fg_{\pm}$ is 
the linear span of the positive (negative) root vectors, implies that 
$\prod_{\alpha\in R_{\od}^+}\alpha\cdot \varphi$ is the restriction of 
an invariant.

\underline{$\fg=\fsl(2|2)$}. In this case the root spaces are 
two-dimensional. A direct calulation proves that $I(\fh^{*})$ 
consists of the polynomials of the form
$$
c+\alpha_{1}\alpha_{2}g+\alpha_{1}^2\alpha_{2}^2\varphi,
$$
where $c\in\Cee$, $g$ is a lilnear combination of the functions of 
the form 
$\frac{\tilde\eps_{1}^n-\tilde\eps_{2}^n}{\tilde\eps_{1}-\tilde\eps_{2}}$ 
for the {\it even} weights $\tilde\eps_{1}$ and $\tilde\eps_{2}$ of the 
identity $\fg$-module and $\varphi\in S(\fh^{*})^W$. 

Proposition 1.4 implies that $\alpha_{1}^2\alpha_{2}^2\varphi$ is the 
restriction of an invariant.  To show that $\alpha_{1}\alpha_{2}g$ is 
the restriction of an invariant, consider the function 
$F(t)=\frac{(t-\tilde\delta_{1})(t-\tilde\delta_{2})}{(t-\tilde\eps_{1})
(t-\tilde\eps_{2})}$, 
where $\tilde\delta_{1}$ and $\tilde\delta_{2}$ are {\it odd} weights 
of the identity $\fg$-module. The coefficients of the power series 
expansion in $t$ of $F(t)$ are expressed in terms of the supertraces 
of powers of the identity representation; hence, are restrictions of 
invariants.

Let $F(t)=\mathop{\sum}\limits_{k\geq 0}t^{-k}\mu_{k}$.
It is easy to check that 
$\mu_{k+2}=\alpha_{1}\alpha_{2}\frac{\tilde\eps_{1}^k- 
\tilde\eps_{2}^k}{\tilde\eps_{1}-\tilde\eps_{2}}$.

\underline{$\fpsl(n|n)$, $n>1$}. Let $\fh\subset \fsl(n|n)$ be the 
Cartan subalgebra, $\tilde\fh=\fh/\Cee\cdot z$, where $z$ is an 
element from the center of $\fsl(n|n)$. It is not difficult to verify 
that $I(\tilde\fh^*)$ can be embedded into $I(\fh^*)$ and the image 
coincides with the set of elements from $I(\fb^*)$ invariant under 
translations in the direction of $z$, i.e., $f$ such that $f(h+tz)=f(h)$. 

Let us continue such a polynomial $f$ to an invariant $F$ from 
$S(\fsl(n|n)^{*})^{\fsl(n|n)}$.  Then $F$ also is invariant 
under translations in the direction of $z$; hence, determines an 
element from $S(\fg^{*})^{\fg}$ whose restriction is equal to $f$. To 
establish this, it suffices to verify that the derivative in the 
direciton of $z$ commutes with the restriction homomorphism onto 
$\tilde\fh$. 

Observe that though there is a nondegenerate invariant supersymmetric 
even bilinear form on $\fg$, by Remark 1.3 this form is NOT related 
with any finite dimensional representation of $\fg$. 

\underline{$\fg= \fosp(2|2n-2)$, $n>1$}.  This Lie superalgebra 
possesses a compatible $\Zee$-grading of depth 1; hence, there is a 
one-to-one correspondence between irreducible finite dimensional 
representations of $\fg$ and irreducible finite dimensional 
representations of $\fg_{\ev}$. Therefore, the arguments applyed for 
$\fgl(m|n)$ are applicable here as well.

Namely, same as for $\fgl(m|n)$, define the algebra
$$
\renewcommand{\arraystretch}{1.2}
\begin{array}{l}
J(\fh^{*})=\{f\in\Cee[[\fh^*]]^{W}\mid f\text{ is a finite linear 
combination of the $e^{\lambda}$}\\
\text{with exponents equal to the weights 
of finite dimensional $\fg_{\ev}$-modules}\\
\text{and $D_{h}f\in(\alpha)$ for any $h\in [\fg_{\od}^{-\alpha}, 
\fg_{\od}^{\alpha}]$ and $\alpha\in R_{\od}^+=\bar R_{\od}^+$}.\}
\end{array}
$$
In the same lines as for $\fgl(m|n)$, we prove that any element from 
$J(\fh^{*})$ is a linear combination of supercharacters.

Let us prove now that any element from $I(\fh^{*})$ is an homogeneous 
component of an element from $J(\fh^{*})$.  Indeed, let $P\in 
I(\fh^{*})$; then $P=P(\eps_{1}, \delta_{1}, \dots ,\delta_{n-1})$, 
where $\eps_{1}, \delta_{1}, \dots ,\delta_{n-1}$ are the weights of 
the standard $\fg$-module. Set
$$
f=P\left(\frac{e^{\eps_{1}}-e^{-\eps_{1}}}{2}, 
\frac{e^{\delta_{1}}-e^{-\delta_{1}}}{2}, \dots 
,\frac{e^{\delta_{n-1}}-e^{-\delta_{n-1}}}{2}\right).
$$
Let us check that $f\in J(\fh^{*})$.  Clearly, $f^{W}=f$.  For 
$\alpha=\eps_{1}-\delta_{1}$ the condition $D_{h}f\in(\alpha)$ is 
equivalent (because the restriction of the invariant form onto $\fh$ 
is proportional to $\eps_{1}^2-\sum\delta_{j}^{2}$) to the fact that 
$\pderf{f}{\eps_{1}}+\pderf{f}{\delta_{1}}\in(\eps_{1}-\delta_{1})$.

Since $P\in I(\fh^{*})$, it follows that 
$\pderf{P}{\eps_{1}}+\pderf{P}{\delta_{1}}=(\eps_{1}-\delta_{1})Q$. 
Hence,
$$
\renewcommand{\arraystretch}{1.2}
\begin{array}{l}
\pderf{f}{\eps_{1}}+\pderf{f}{\delta_{1}}= 
\pderf{P}{\eps_{1}}\cosh(\eps_{1})+\pderf{P}{\delta_{1}}
\cosh(\delta_{1})=\\ \pderf{P}{\eps_{1}}\cosh(\eps_{1})+ 
(\sh(\eps_{1})-\sinh(\delta_{1}))Q\cosh(\delta_{1})- 
\pderf{P}{\eps_{1}}\cosh(\delta_{1})\\
\pderf{P}{\eps_{1}}(\cosh(\eps_{1})-\cosh(\delta_{1}))+ 
(\sinh(\eps_{1})-\sinh(\delta_{1}))Q\cosh(\delta_{1})\in (\alpha).
\end{array}
$$
Therefore, $f\in J(\fh^{*})$ and its homogeneous component of the 
least degree --- equal to $P$ --- is the restriction of an invariant.

\underline{$\fg= \fosp(2m+1|2n)$}.  For a basis of $\fh^{*}$ we take 
the weights $\eps_{1}, \dots , \eps_{m}; \delta_{1}, \dots 
,\delta_{n}$ of the identity representation (of the two weights 
$\pm\eps_{i}$ and $\pm\delta_{j}$ we select one).  Observe that the 
odd roots $\pm\delta_{j}$ of $\fg$ are collinear to the even ones. 

The Weyl group separately permutes the $\eps$'s and the $\delta$'s and 
changes their signs.  The restriction of the invariant form onto $\fh$ 
is proportional to $\sum\eps_{i}^2-\sum\delta_{j}^{2}$; hence, for 
$\alpha= \eps_{i}-\delta_{j}$ (observe that $\alpha\in \tilde R_{\od}$) 
the condition $D_{h}f\in(\alpha)$ is equivalent to the fact that 
$\pderf{f}{\eps_{i}}+\pderf{f}{\delta_{j}}\in(\eps_{i}-\delta_{j})$, 
which, in turn, means that $f$ does not depend on $t$ after 
substitution $\eps_{i}=\delta_{j}=t$. Clearly,
$$
f=f(\eps_{1}^2, \dots , \eps_{m}^2, \delta_{1}^2, \dots 
,\delta_{n}^2)
$$
and, therefore, is a supersymmetric polynomial in the sence of 
\cite{P} and as such can be expressed via the coefficients of the 
rational function
$$
F(t)=\frac{\prod(t^2-\delta_{j}^2)}{\prod(t^2-\eps_{i}^2)}.
$$ 
These coefficients are expressed via the sums 
$\mathop{\sum}\limits_{i}\eps_{i}^k-\mathop{\sum}\limits_{j}\delta_{j}^{k}$
--- the powers of the supertrace the identity representation.

\underline{$\fg= \fosp(2m|2n)$}.  For a basis of $\fh^{*}$ we take the 
same basis $\eps_{1}, \dots , \eps_{m}$; $\delta_{1}, \dots 
,\delta_{n}$ as in the preceding case.  The restriction of any 
$W$-invariant polynomial is of the form 
$$
f=P(\eps_{1}^2, \dots , \eps_{m}^2, \delta_{1}^2, \dots 
,\delta_{n}^2)+\eps_{1}\cdot \dots \cdot \eps_{m}\cdot Q(\eps_{1}^2, 
\dots , \eps_{m}^2, \delta_{1}^2, \dots ,\delta_{n}^2).
$$
Moreover, the restriction of any $W$-invariant polynomial does not 
depend on $t$ after substitution $\eps_{i}=\delta_{j}=t$, whereas 
after such a substitution $P$ is of an even degree wrt $t$ and the 
second summand is of an odd degree.  This means that both summands do 
not depend on $t$; hence, $P\in I(\fh^{*})$ and $\eps_{1}\dots 
\eps_{m}Q\in I(\fh^{*})$ and 
$$
Q=\prod (\eps_{i}^2-\eps_{i}^2)Q_{1}\; \text{ for some polynomial } 
Q_{1}\eqno{(*)}.
$$

The same arguments as in the preceding case show that $P$ is the 
restriction of an invariant.  Let us prove that $\eps_{1}\dots 
\eps_{m}Q$ is also the restriction of an invariant.  To this end, 
select in $\fh$ the right dual basis $e_{1}, \dots , e_{m}$; $f_{1}, 
\dots ,f_{n}$ to $\eps_{1}, \dots , \eps_{m}$; $\delta_{1}, \dots 
,\delta_{n}$.  For the system of simple roots take
$$
\delta_{1}-\delta_{2}, \dots , \delta_{n-1}-\delta_{n}; 
\delta_{n}-\eps_{1}, \eps_{1}-\eps_{2}, \dots , \eps_{m-1}-\eps_{m}, 
\eps_{m-1}+\eps_{m}.
$$
Let $\Lambda\in \fh^{*}$ be such that $\Lambda(e_{i})=\lambda_{i}$ and 
$\Lambda(f_{j})=\mu_{j}$. Then due to \cite{K1} for the 
representation with highest wight $\Lambda$ to be a finite dimensional 
one, the coordinates of the highsest weight should satisfy the 
following conditions:
$$
\renewcommand{\arraystretch}{1.2}
\begin{array}{l}
\mu_{j}\in\Zee\text{ and }\mu_{1}\geq \mu_{2}\geq \dots \geq 
\mu_{n}\geq m;\\
\lambda_{1}-\lambda_{2}\in\Zee_{+},\; \dots ,\; 
\lambda_{m-1}-\lambda_{m}\in\Zee_{+},\; 
\lambda_{m-1}+\lambda_{m}\in\Zee_{+}.
\end{array}
$$
If $\Lambda$ satisfies $\lambda_{m}\neq 0$, then $\Lambda$ is a 
typical weight. Indeed, by Kac, \cite{K1}, the typicality condition is
$$
(\Lambda+\rho)(h_{\alpha})\neq 0\text{ for any }\alpha\in R_{\od}^+
$$
or, equivalently,
$$
(\Lambda+\rho)(e_{i})\neq 0\text{ for any }i\text{ and }
(\Lambda+\rho)(f_{j})\neq 0\text{ for any }j.
$$
As is not difficult to verify, 
$$
\rho_{0}=(m-1)\eps_{1}+ \dots +\eps_{m-1}+n\delta_{1}+ \dots 
+2\delta_{n};\; \; \rho_{1}=m\sum\delta_{i};
$$ 
hence, 
$$
(\Lambda+\rho)(e_{i})=\lambda _{i}+m-i\neq 0\text{ and }
(\Lambda+\rho)(f_{j})=\mu _{j}-m+m-j+1\neq 0
$$
as was required for typicality.

Let $\Lambda$ be a typical weight.  Then by \cite{K2} the 
supercharacter of the irreducible module $L^{\Lambda}$ with the 
highest weight $\Lambda$ is
$$
\renewcommand{\arraystretch}{1.4}
\begin{array}{l}
\sch L^{\Lambda}=\frac{\mathop{\prod}\limits_{\alpha\in 
R_{\od}^+}\left(e^{\alpha/2}- e^{-\alpha/2}\right) } 
{\mathop{\prod}\limits_{\alpha\in 
R_{\ev}^+}\left(e^{\alpha/2}-e^{-\alpha/2}\right)} 
\sum\eps(w)e^{w(\Lambda+\rho)}=\\
\mathop{\prod}\limits_{\alpha\in R_{\od}^+}\left(e^{\alpha/2}- 
e^{-\alpha/2}\right) \mathop{\prod}\limits_{\alpha\in 
R_{\ev}^+}\left(e^{\alpha/2}-e^{-\alpha/2}\right)^{-1} 
\sum\eps(w)e^{w(\Lambda+\rho_{0}-\rho_{1})}=\\
\mathop{\prod}\limits_{\alpha\in R_{\od}^+}\left(e^{\alpha/2}- 
e^{-\alpha/2}\right) \ch L_{0}^{\Lambda-\rho_{1}},
\end{array}
$$
where $L_{0}^{\Lambda-\rho_{1}}$ is the irreducible $\fg_{\ev}$-module 
with the highest weight $\Lambda-\rho_{1}$. 

Consider $\Lambda$ such that $\lambda_{m}\neq 0$ and $\mu_{n}\geq 
m$, i.e., consider the highest weights of the typical irreducible 
finite dimensional modules. Then $\Lambda-\rho_{1}$ runs the highsest 
weights of $\fg_{\ev}$ for which $\lambda_{m}\neq 0$. Let us demand 
that $\lambda_{m}$ were half-integer; then $\lambda_{m}\neq 0$ 
holds automatically. Let $T\in S(\fh^{*}_{\ev})^{W}$ be an invariant 
polynomial of the form
$$
T=\eps_{1}\cdot \dots \cdot \eps_{m}\cdot Q_{1}(\eps_{1}^2, 
\dots , \eps_{m}^2, \delta_{1}^2, \dots ,\delta_{n}^2)\; \text{ for 
$Q_{1}$ defined in $(*)$}.
$$
Consider
$$
\renewcommand{\arraystretch}{1.4}
\begin{array}{l}
\tilde T=T(e^{\eps_{1}/2}-e^{-\eps_{1}/2}, \dots , 
e^{\eps_{m}/2}-e^{-\eps_{m}/2}; e^{\delta_{1}/2}-e^{-\delta_{1}/2}, 
\dots , e^{\delta_{n}/2}-e^{-\delta_{n}/2})=\\
\mathop{\prod}\limits_{i=1}^m(e^{\eps_{i}/2}-e^{-\eps_{i}/2}) \cdot 
S(2\sinh(\eps_{1}/2), \dots , 2\sinh(\eps_{m}/2); 2\sinh(\delta_{1}/2) , 
\dots , 2\sinh(\delta_{n}/2)).
\end{array}
$$
All the weights of the first factor are half-integer and linealy 
independent, hence, it is a linear combination of the $e^\chi$, where 
the coordinates $\chi_{i}$ are half-integers and nonzero (actually, 
this $\chi$ is the character of a spinor-like representation).

At the same time all the weights of the second factor are integer; so 
$T$ is a linear combination of the $e^\chi$ with $\chi_{m}\neq 0$ and 
half-integer.  This expression is a linear combination of characters 
of finite dimensional representations
$$
\tilde T=\sum C_{\chi}\ch L_{0}^\chi
$$ 
of $\fg_{\ev}$. Since all the weights $\mu$ that 
contribute to $S$ are half-integer with $\mu_{m}\neq 0$, then 
all the coordinates of $\chi$ are also half-integer and $\chi_{m}\neq 0$.

By multiplying both sides of the last equation by 
$L=\mathop{\prod}\limits_{\alpha\in R_{\od}^+}\left(e^{\alpha/2}-
e^{-\alpha/2}\right)$ we obtain
$$
L\cdot \tilde T=\sum C_{\chi}L\cdot \ch L_{0}^\chi=C_{\chi}\sch 
L^{\chi+\rho_{1}}.
$$
Therefore, the lowest component of $L\tilde S$ is the restriction of an 
invariant. But this lowest component is $\eps_{1}\cdot \dots \cdot \eps_{m}
\cdot Q(\eps_{1}^2, 
\dots , \eps_{m}^2, \delta_{1}^2, \dots ,\delta_{n}^2)$.

\underline{$\fg= \fosp_{\alpha}(4|2)$}.  In this case 
$\fg_{\ev}=\fsl_{1}(2)\oplus\fsl_{2}(2)\oplus\fsl_{3}(2)$ (the sum of 
three copies of $\fsl(2)$, numbered to distinguish them) and 
$\fg_{\od}=V_{1}\otimes V_{2}\otimes V_{3}$, where $V$ is the identity 
representation. Let $\pm\eps_{1}$, $\pm\eps_{2}$, $\pm\eps_{3}$ be the 
weights of each of the components of $\fg_{\ev}$ in its identity module. 
Then the root system of $\fg$ is as follows
$$
R_{\ev}=\{\pm 2\eps_{1}, \; \pm 2\eps_{2}, \; \pm2\eps_{3}\},\; \; 
R_{\od}=\{\pm \eps_{1} \pm \eps_{2}\pm \eps_{3}\}.
$$
The restriction of the nondegenerate invariant supersymmetric even 
bilinear form onto the Cartan subalgebra is of the form
$$
\frac{1}{\lambda_{1}}\eps_{1}^{2}+ \frac{1}{\lambda_{2}}\eps_{2}^2+
\frac{1}{\lambda_{3}}\eps_{3}^{2}, \text{ where 
}\lambda_{1}=-(1+\alpha), \; \lambda_{2}=1, \; \lambda_{3}=\alpha.
$$
From the condition $(h_{\alpha}, h)=\alpha(h)$ for $h\in \fh$ we 
deduce that if $\alpha =\theta_{1}\eps_{1}+\theta_{2}\eps_{2}+
\theta_{3}\eps_{3}$, where $\theta_{i}=\pm 1$, then 
$$
h_{\alpha}=\theta_{1}\lambda_{1}H_{1}+\theta_{2}\lambda_{2}H_{2}+
\theta_{3}\lambda_{3}H_{3}
$$
for the dual basis of $\fh$, i.e., $\eps_{i}(H_{j})=\delta_{ij}$.

A direct calculation shows that 
$$
D_{h_{\alpha}}f=\theta_{1}\lambda_{1}\pderf{f}{\eps_{1}}+
\theta_{2}\lambda_{2}\pderf{f}{\eps_{2}}+
\theta_{3}\lambda_{3}\pderf{f}{\eps_{3}}.
$$

The Weyl group is isomorphic to $(\Zee/2)^{3}$; it acts by changings 
the signs of the $\eps_{i}$.

Let us describe now the algebra $I(\fh^*)$. An easy calculation shows 
that $I(\fh^*)$ consists of the elements of the form
$$
f=\psi(\frac{1}{\lambda_{1}}\eps_{1}^{2}+ \frac{1}{\lambda_{2}}\eps_{2}^2+
\frac{1}{\lambda_{3}}\eps_{3}^{2})+\prod_{\alpha\in R_{\od}^+}
\alpha\cdot \varphi,
$$
where $\varphi\in S(\fh^*)^W$.  Clearly, the first summand is the 
restriction of an invariant. Let us prove that the second summand is 
also the restriction of an invariant. 

For a system of simple roots take $\{\eps_{1}+\eps_{2}+\eps_{3},\; 
-2\eps_{2}, \; -2\eps_{3}\}$.  Applying Proposition 1.5 to $\fg$ we 
see that any element of the form $\mathop{\prod}\limits_{\alpha\in 
R_{\od}^+} \alpha^2\cdot \varphi$, where $\varphi\in S(\fh)^W$ is 
the restriction of an invariant.  Since $\fg$ possesses a 
nondegenerate invariant supersymmetric even bilinear form and the 
Harish--Chandra homomorphism applies, any element from $S(\fh)^W$ 
can be expressed in the form
$$
\frac{f}{Q^2},\text{ where $f$ is the image of an element from the 
center of $U(\fg)$ and $Q=\prod_{\alpha\in R_{\od}^+}h_{\alpha}$}.
$$
This implies that the irreducible finite dimensional $\fg$-module 
with highest weight $\Lambda$ is typical if 
$$
\prod_{\alpha\in R_{\od}^+}\Lambda(h_{\alpha})\neq 0.
$$

\begin{Remark} In \cite{K2} this statement is conjectured, while the 
conditions sufficient for typicality of the modules over 
$\fosp_{\alpha}(4|2)$, $\fag_{2}$ and $\fab_{3}$ (these Lie 
superalgebras are differently baptized there) given there are 
faulty.)\end{Remark}

As for $\fosp(2n|2m)$, the formula for the supercharacter of the 
typical module can be expressed in the form
$$
\sch L^\Lambda=\prod_{\alpha\in R_{\od}^+}\left(e^{\alpha/2}-
e^{-\alpha/2}\right)\ch L_{0}^{\Lambda-\rho_{1}}.
$$
In the chosen system of simple roots
$$
\rho_{0}=\eps_{1}-\eps_{2}-\eps_{3};\; \; \rho_{1}=2\eps_{1}.
$$
If $\Lambda= \chi_{1}\eps_{1}+\chi_{2}\eps_{2}+\chi_{3}\eps_{3}$ is 
the highest weight, then $L^\Lambda$ is of finite dimension whenever 
$$
\chi_{1}\geq 2, \; \chi_{1}\in\Zee_{+};\; \;  \chi_{2}, \chi_{3}\in\Zee_{-};
$$
the module $L^\Lambda$ is typical if additionally
$$
\prod_{\alpha\in R_{\od}^+}\Lambda(h_{\alpha})\neq 0.
$$

Fix $\chi_{2}$ and $\chi_{3}$.  The equation for $\chi_{1}$ obtained 
has finitely many solutions.  Hence, by selecting $\chi_{1}$ 
sufficiently large we deduce that for arbitrary $\chi_{2}$ and 
$\chi_{3}$ (such that $\chi_{2}, \chi_{3}\in\Zee_{-}$) the module 
$L^\Lambda$ with $\Lambda= 
\chi_{1}\eps_{1}+\chi_{2}\eps_{2}+\chi_{3}\eps_{3}$ is finite 
dimensional and typical.

Now, represent $\fg_{\ev}$ in the form 
$\fg_{\ev}=\fsl(2)\oplus\tilde\fg$ and set $\fh=\Span(H_{1})\oplus 
\tilde\fh$, respectively.  According to Bernstein \cite{B}, any element 
$P$ from $S(\tilde\fh^*)^{\Zee/2\times\Zee/2}$ is the component of the 
least degree of the linear combination $T$ of characters of finite 
dimensional representations of $\tilde\fg$.  By Proposition 1.6, 
$\eps_{1}^k$ is a linear combination $\tilde T$ of characters of finite 
dimensional representations of $\fsl(2)$ whose highest weight is 
sufficiently big.  Then the product $\eps_{1}^k\cdot P$ is a linear 
combination of characters of finite dimensional 
representations of $\fg$ that enter $\tilde T\cdot T$ and each of which 
satisfies requirements for finite dimension and typicality.

Therefore, $\prod_{\alpha\in R_{\od}^+}\left(e^{\alpha/2}-
e^{-\alpha/2}\right)\eps_{1}^k\cdot P$ is a linear 
combination of supercharacters on $\fg$ and its lowest component is 
equal to $\prod_{\alpha\in R_{\od}^+}\alpha\cdot\eps_{1}^k\cdot P$.

\underline{$\fg= \fag_{2}$}.  In this case 
$\fg_{\ev}=\fg_{2}\oplus\fsl(2)$ and $\fg_{\od}=R(\pi_{1})\otimes V$, 
where $R(\pi_{1})$ is the first fundamental representation of 
$\fg_{2}$ (see \cite{OV} or \cite{Bu1}) and $V$ is the identity 
representation of $\fsl(2)$.  We select a basis $H_{1}$, $H_{2}$, 
$H_{3}$, in the Cartan subalgebra of $\fg_{2}$ so that 
$H_{1}+H_{2}+H_{3}=0$; let the $\lambda_{i}$ be the linear forms such 
that $\lambda_{i}(H_j)=-1$ if $i\neq j$ and $\lambda_{i}(H_i)=2$.

Let $\pm\delta$ be the weights of the identity $\fsl(2)$-module; 
select the basis element $H$ of the Cartan subalgebra of $\fsl(2)$ so 
that $\delta(H)=1$.  Then the root system of $\fg$ is as follows
$$
R_{\ev}=\{\lambda_{i}-\lambda_{j}; \; \pm \lambda_{i}; \; 
\pm2\delta\},\; \; R_{\od}=\{\pm\lambda_{i}\pm \delta; \; \pm\delta\}.
$$
For the system of simple roots select
$$
\lambda_{1}+\delta,\; \lambda_{2}; \; \lambda_{3}-\lambda_{2}.
$$
Then 
$$
\renewcommand{\arraystretch}{1.2}
\begin{array}{l}
R_{\ev}^+=\{\lambda_{2}, \, \lambda_{3}, \, -\lambda_{1}; \;
\lambda_{2}-\lambda_{1}, \; 
\lambda_{3}-\lambda_{2}, \; \lambda_{3}-\lambda_{1},\; 
2\delta\},\\
R_{\od}^+=\{\lambda_{1}+\delta,\; \lambda_{2}+\delta,\; 
\lambda_{3}+\delta,\; -\lambda_{1}+\delta,\; -\lambda_{2}+\delta,\; 
-\lambda_{3}+\delta,\; \delta\}.
\end{array}
$$

Observe that $\delta\in R_{\od}^+$ while $2\delta\in R_{\ev}^+$; 
hence, $\tilde R_{\od}= R_{\od}^+\setminus\{\delta\}$.  The Weyl group 
is isomorphic to $(S_{3}\circ\Zee/2)\times \Zee/2$, where $S_{3}$ 
permutes the $\lambda_i$, the first $\Zee/2$ simultaneusly changes the 
signs of all the $\lambda_i$, the second $\Zee/2$ changes the sign of 
$\delta$.  It is not difficult to verify that
$$
S(\fh^*)^W=\Cee[\lambda_{1}^2+\lambda_{2}^2+\lambda_{3}^2, \; \; 
(\lambda_{1}\lambda_{2}\lambda_{3})^2, \; \; \delta^2].
$$
The restriction of the nondegenerate invariant supersymmetric even 
bilinear form $\fh$ is proportional to
$$
3\delta^2-2(\lambda_{1}^2+\lambda_{2}^2+\lambda_{3}^2).
$$

From the condition $(h_\alpha, h)=\alpha(h)$ for $h\in \fh$ 
setting $\alpha=\lambda_{3}+\delta$ and selecting $H_1$, $H_2$, 
$H$ for a basis  of $\fh$ we deduce that
$$
h_\alpha=-H_1-H_2-2H.
$$
Therefore, the condition $D_{h_\alpha}f\in (\alpha)$ is equivalent to 
the following one:
$$
\pderf{f}{\lambda_{1}}+\pderf{f}{\lambda_{2}}+2\pderf{f}{\delta}+
\in(\delta-\lambda_{1}-\lambda_{2}).
$$
The other conditions of a similar type follow from $W$-invariance. 
Direct calculations demonstrate that $I(\fh^{* })$ consists of the 
elements of the form
$$
f=\psi(3\delta^2-2(\lambda_{1}^2+\lambda_{2}^2+\lambda_{3}^2))+
\prod_{\alpha\in \tilde R_{\od}^+}
\alpha\cdot \varphi,
$$
where $\varphi\in S(\fh^*)^W$.  

The same arguments as in the study of $\fosp_\alpha(4|2)$ proove with 
the help of Proposition 1.5 that if $\prod_{\alpha\in 
R_{\od}^+}\Lambda(h_{\alpha})\neq 0$ then the module $L^\Lambda$ is a 
typical one.  In this case the formula for the supercharacter can be 
expressed in the form (\cite{K2})
$$
\sch L^\Lambda=\prod_{\alpha\in R_{\od}^+}\left(e^{\alpha/2}-
e^{-\alpha/2}\right)\prod_{\alpha\in R_{\ev}^+}\left(e^{\alpha/2}-
e^{-\alpha/2}\right)^{-1}\mathop{\sum}\limits_{w\in W}\eps'(w)
e^{w(\Lambda+\rho)},
$$
where $\eps'(w)=(-1)^{N(w)}$ and $N(w)$ is the parity of the number of 
reflections in all even roots except $2\delta$.

Let $W_{1}=S_{2}\circ\Zee/2$ be the Weyl group of $\fg_{2}$; let 
$\rho_{0}{}'$ be the halfsum of the positive roots for $\fg_{2}$ and  
$\rho_{0}{}''=\delta$ the halfsum of the positive roots for $\fsl(2)$; 
let $\rho_{1}=\frac72\delta$ be the halfsum of the positive odd roots.

If $w$ is the reflection in $2\delta$, then $\eps'(w)=1$. Therefore 
the formula for the supercharacter can be represented in the 
following form
$$
\renewcommand{\arraystretch}{1.4}
\begin{array}{l}
\sch L^\Lambda=\mathop{\prod}\limits_{\alpha\in \tilde 
R_{\od}^+}\left(e^{\alpha/2}- 
e^{-\alpha/2}\right)\mathop{\prod}\limits_{\alpha\in \tilde 
R_{\ev}^+}\left(e^{\alpha/2}- 
e^{-\alpha/2}\right)^{-1}\left(e^{\delta/2}- 
e^{-\delta/2}\right)\left(e^{\delta}- e^{-\delta}\right)^{-1}\times\\
\Big(\mathop{\sum}\limits_{w\in 
W_1}\eps(w)e^{w(\Lambda'+\rho'{}_{0})}\Big)\cdot 
\Big(e^{\Lambda''+\rho''{}_{0}-\rho_{1}}+ 
e^{-(\Lambda''+\rho''{}_{0}-\rho_{1})}\Big)=\\
\mathop{\prod}\limits_{\alpha\in \tilde R_{\od}^+}\left(e^{\alpha/2}-
e^{-\alpha/2}\right)\left(e^{\delta/2}+
e^{-\delta/2}\right)^{-1}\ch L_{0}^{\Lambda{}'}
\Big(e^{\Lambda''+\rho''{}_{0}-\rho_{1}}+
e^{-(\Lambda''+\rho''{}_{0}-\rho_{1})}\Big),
\end{array}
$$
where $\Lambda{}'$ is the restriction of $\Lambda$ onto the Cartan 
subalgebra of $\fg_{2}$ and $\Lambda{}''$ is the same for $\fsl(2)$ 
under the proviso: 
$$
\Lambda{}''(H)\geq 7.
$$
Under this condition $\dim L^\Lambda<\infty$. To see this, it saffices 
to take the $\fg_{\ev}$-module $\ch L_{0}^{\Lambda}$ with the same 
highest weight and the highest weight vector $v$; then in the induced 
module $\Ind_{\fg_{\ev}}^{\fg}L_{0}^{\Lambda}$, consider the 
submodule generated by $\mathop{\prod}\limits_{\alpha\in  
R_{\od}^+}h_\alpha\cdot v$. (This submodule is currently called 
sometimes in the literature ``Kac' module''.)

Let $P\in I(\fh^{*})$. Then, as we have shown, 
$$
P=P_{1}(3\delta^2-2(\lambda_{1}^2+\lambda_{2}^2+\lambda_{3}^2))+
\mathop{\prod}\limits_{\alpha\in R_{\od}^+}\alpha\cdot P_{2},
$$
where $P_{2}\in S(\fh^*)^W$.  As we have already mentioned, 
$3\delta^2-2(\lambda_{1}^2+\lambda_{2}^2+\lambda_{3}^2$ is the 
restriction of $\str~\ad~X$.  Let us show that 
$\mathop{\prod}\limits_{\alpha\in R_{\od}^+}\alpha\cdot P_{2}$ is also 
the restriction of an invariant.  Since $P_{2}\in S(\fh^*)^W$, we may 
assume that $P_{2}=Q_{2}\cdot \delta^{2k}$, where $Q_{2}\in S\left 
((\fg_{2}\cap \fh)^*\right)^{W_{1}}$.

According to Bernstein (\cite{B}) $Q_{2}$ is the lowest component of a 
linear combination $T$ of characters of irreducible $\fg_{2}$-modules.

By Remark 1.6, $\delta^{2k}$ is the lowest component of a linear 
combination $T'$ of functions of the form $\frac{\cosh 
(2k-5)\delta}{\cosh \delta}$.  So $P$ is the lowest component of the 
linear combination $TT'$.  If we now fix $\Lambda{}'$, then the 
equation for $\Lambda{}''$ obtained from typicality conditions has 
finitely many solutions.  Hence, by selecting sufficiently large 
numbers $k$ in $T'$ we may assume that any highest weight of the 
linear combination $TT'$ satisfies the conditions for finite dimension 
and typicality.

This implies that the lowest component of 
$\mathop{\prod}\limits_{\alpha\in \tilde R_{\od}^+}\left(e^{\alpha/2}- 
e^{-\alpha/2}\right)\cdot TT'$, equal to $P_{2}$, is the restriction 
of an invariant.

\underline{$\fg=\fab_{3}$}.  In this case 
$\fg_{\ev}=\fo(7)\oplus\fsl(2)$ and $\fg_{\od}=\spin(7)\otimes \id$.  
Let $\pm\eps_{1}$, $\pm\eps_{2}$, $\pm\eps_{3}$ be the weights of the 
standard $\fo(7)$-module and $\pm\frac12\delta$ the weights of the 
identity representation of $\fsl(2)$.  Then
$$
R_{\ev}=\{\pm\eps_{i}; \; \pm\eps_{i}\pm\eps_{j}; \; \pm\delta\},\; \; 
R_{\od}=\{\frac12(\pm\eps_{1}\pm\eps_{2}\pm\eps_{3}\pm\delta)\}.
$$
For the system of simple roots select
$$
\frac12(\delta-\eps_{1}-\eps_{2}-\eps_{3}), \; \eps_{1}-\eps_{2}, \; 
\eps_{2}-\eps_{3}, \; \eps_{3}.
$$
Then 
$$
\renewcommand{\arraystretch}{1.2}
\begin{array}{l}
R_{\ev}^+=\{\eps_{2}; \; \eps_{2}, \; 
\eps_{3}, \; \eps_{i}\pm\eps_{j}\text{ for }i<j,\; 
\delta\},\\
R_{\od}^+=\{\frac12(\delta\pm\eps_{1}\pm\eps_{2}\pm\eps_{3})\}.
\end{array}
$$
Clearly, $\rho_{1}=2\delta$.

The restriction of the nondegenerate invariant supersymmetric even 
bilinear form $\str(\ad~x\cdot \ad~y)$ onto the Cartan subalgebra is 
proportional to
$$
3(\eps_{1}^2+\eps_{2}^2+\eps_{3}^2)-\delta^2.
$$

From the condition $(h_\alpha, h)=\alpha(h)$ for $h\in \fh$ we deduce 
for $\alpha=\frac12(\eps_{1}+\eps_{2}+\eps_{3}+\delta)$ that
$$
h_\alpha=\frac16(H_1+H_2+H_3)-\frac12H,
$$
where $H_1$, $H_2$, $H_3$ and $H$ is the basis of $\fh$ dual to 
$\eps_{1}$, $\eps_{2}$, $\eps_{3}$ and $\delta$.
Therefore, the condition $D_{h_\alpha}f\in (\alpha)$ is equivalent to 
the following one:
$$
\pderf{f}{\eps_{1}}+\pderf{f}{\eps_{2}}+\pderf{f}{\eps_{3}}- 
3\pderf{f}{\delta} \in(\delta+\eps_{1}+\eps_{2}+\eps_{3}).
$$
The other conditions of a similar type follow from $W$-invariance.  
Observe that $W$ is isomorphic to $(S_{3}\circ(\Zee/2)^3)\times 
\Zee/2$.  Direct calculations demonstrate that $I(\fh^{* })$ consists 
of the elements of the form
$$
P_{1}(\mu_{-2}, \mu_{2})+ \mathop{\prod}\limits_{\alpha\in 
R_{\od}^+}\alpha\cdot P_{2},
$$
where the $\mu_i$ are the coefficients of $t^i$ in the power series 
expansion of the rational function 
$\frac{\mathop{\prod}\limits_{\alpha\in 
R_{\od}}(t-\alpha)}{\mathop{\prod}\limits_{\alpha\in 
R_{\ev}}(t-\alpha)}$ and $P_{2}\in S(\fh^*)^W$.

It is clear that $\mu_{\pm 2}$  are the restrictions of invariants; we 
prove that $\mathop{\prod}\limits_{\alpha\in R_{\od}^+}\alpha\cdot 
P_{2}$ is the restriction of an invariant by the same arguments as in 
the above two last cases.

\underline{$\fg=\fvect(0|n)$, $\fsvect(0|n)$ or 
$\widetilde{\fsvect}(0|2n)$}. Let 
$\fg=\mathop{\oplus}\limits_{i\geq -1}\fg_{i}$ be the standard 
$\Zee$-grading; $\eps_{1}$, \dots , $\eps_{n}$ the weights of the 
$\fg_{0}$ -module $\fg_{-1}$. Then it is not difficult to verify that 
$[\fg_{\od}^\alpha,  \fg_{\od}^{-\alpha}]=\ker \alpha$ for 
$\alpha=\eps _{i}$ and the condition $D_{h_\alpha}f\in (\alpha)$ 
implies that $f|_{\ker \alpha}= const=c$. 

Hence, $f-c$ is divisible by $\alpha$ and by $W$-symmetry
it is divisible by $\eps_{1}\dots \eps_{n}$ since the factors are 
mutually prime. Therefore, any element from $I(\fh^{*})$ is of the form
$$
c+\eps_{1}\dots \eps_{n}\cdot P, \text{ where }P\in S(\fh^*)^W.
$$
Proposition 1.4 implies that such an element is the restriction of 
an invariant.

\underline{$\fg=\fpe(n)$}.  Since $\fg$ possesses a compatible 
$\Zee$-grading of depth 1, there is a one-to-one correspondence 
between the irreducible $\fg$-modules and the irreducible 
$\fg_{\ev}$-modules. Therefore, arguments similar to the ones applied 
for $\fg=\fgl$ show that any of the elements from $J(\fh^{*})$ (which 
is similarly defined) are linear combinations of supercharacters of 
finite dimensional representations. Now, let $P\in I(\fh^{*})$. If 
$\pm\eps_{1}$, \dots , $\pm\eps_{n}$ are the weights of the identity
$\fg$-module, then
$$
R_{\od}=\{\pm(\eps_{i}+\eps_{j})\text{ for }i\neq j, \; \; 
-2\eps_{i}\}.
$$
Let $e_{1}$, \dots , $e_{n}$ be the basis of $\fh$ dual to $\eps_{1}$, 
\dots , $\eps_{n}$.  Then for $\alpha=\eps_{i}+\eps_{j}$ for $i\neq j$ 
and $h_\alpha=e_{i}-e_{j}$ the condition $D_{h_\alpha}f\in (\alpha)$ 
means that
$$
\pderf{P}{\eps_{i}}-\pderf{P}{\eps_{j}}\in(\eps_{i}+\eps_{j}).
$$
Equivalently, one can say that $P$ does not depend on $t$ after the 
substitution $\eps_{i}=-\eps_{j}=t$.

The same arguments as for $\fg=\fgl(m|n)$ prove that if $P\in 
I(\fh^{*})$, then 
$$
f=P(e^{\eps_{1}/2}- e^{-\eps_{1}/2}, \dots , e^{\eps_{n}/2}- 
e^{-\eps_{n}/2})\in J(\fh^{*})
$$ 
and the homogeneous component of $f$ of the least degree is equal to 
$P$.

\underline{$\fg=\fs\fpe(n)$, $n\neq 4$}.  The answer is the same for 
all simple Lie superalgebras $\fs\fpe(n)$, but the proof is different, 
as we will see. The exceptional case $n=4$ is considered separately.

Let us show that the 
restriction map $I(\fh^{*})\tto I(\tilde \fh^{*})$, where $\fh$ is the
Cartan subalgebra in $\fpe(n)$, while $\tilde \fh$ is same in $\fs\fpe(n)$, 
is surjective.

In notations of the above case, observe that
$$
I(\fh^{*})=\Cee[\Delta_{1}, \Delta_{3}, \dots , \Delta_{2k+1}, 
\dots ],
$$
where $\Delta_{l}=\mathop{\sum}\limits_{i=1}^n\eps_{i}^l$.  Let 
$\tilde\eps_{i}$ be the restriction of $\eps_{i}$ onto $\tilde \fh$.  
Then
$$
I(\tilde \fh^{*})=\Cee[\sigma_{1}, \sigma_{3}, \dots , \sigma_{2k+1}, 
\dots ],
$$
where $\sigma_{l}=\mathop{\sum}\limits_{i=1}^{n-1}\tilde\eps_{i}^l$.

Let $f\in I(\tilde\fh^{*})$, then $f=F(\sigma_{1}, \sigma_{3}, \dots , 
\sigma_{2k+1}, 
\dots )$. Consider the difference
$$
f-F(\sigma_{1}-\sigma_{1}, \sigma_{3}-\sigma_{1}^3, \dots , 
\sigma_{2k+1}-\sigma_{1}^{2k+1}, 
\dots )
$$
Under substitution $\sigma_{1}=0$ the difference vanishes, hence, is 
divisible by $\sigma_{1}$. Since $\sigma_{2k+1}-\sigma_{1}^{2k+1}$ is 
the image of an element from $I(\fh^{*})$ (namely, of 
$\Delta_{2k+1}$), we may assume that $f$ is divisible by $\sigma_{1}$, 
i.e., $f=\sigma_{1}g$.

Let $\alpha=\eps_{n-1}+\eps_{n}$.  The condition $D_{h_\alpha}f\in 
(\alpha)$ means that the restriction of $f$ onto $\ker \alpha$ is 
invariant with respect to translations by $h_\alpha$.  Since 
$\sigma_{1}(h_\alpha)\neq 0$, any element from  $\ker \alpha$ can be 
expressed in the form 
$$
h+th_\alpha,\text{ where }h\in \ker \alpha\cap \ker \sigma_{1}.
$$
Therefore,
$$
f(h+th_\alpha)=\sigma_{1}(h+th_\alpha)g(h+th_\alpha)=
t\sigma_{1}(h_\alpha)g(h+th_\alpha)=f(h),
$$
hence, $f(h+th_\alpha)=0$.  So, $f$ is divisible by $\alpha$ and 
$W$-symmetry and the fact that for $n\neq 4$ the linear functions 
$\alpha$ from $\tilde R_{\od}^+$ are mutually prime we deduce that $f$ 
is divisible by $\mathop{\prod}\limits_{\alpha\in\tilde 
R_{\od}}\alpha$, i.e.,
$$
f=\mathop{\prod}\limits_{\alpha\in\tilde R_{\od}}\alpha\cdot \varphi, 
\text{ where } \varphi\in S(\tilde\fh^*)^W.
$$ 
By Proposition 1.4 such an element is the restriction of an invariant.

\underline{$\fg=\fs\fpe(4)$}.  In this subcase the root spaces are 
2-dimensional.

Let $\tilde\eps_{i}$ be the restriction of $\eps_{i}\in\fh^*$ onto 
$\tilde\fh$, where $\fh$ is the Cartan subalgebra in $\fpe(4)$ and 
$\tilde\fh$ is same in $\fspe(4)$.  Then 
$\alpha_i=\tilde\eps_{1}+\tilde\eps_{2}+\tilde\eps_{3}-\tilde\eps_{i}$. 
 The direct calculations show that $I(\tilde\fh^{*})$ consists of 
the elements of the form
$$
c+\alpha_{1}\alpha_{2}\alpha_{3}g+
(\alpha_{1}\alpha_{2}\alpha_{3})^2\varphi,
$$
where $c\in\Cee$, $g$ is any linear combination of the coefficients of 
the rational function
$$
F(t)=\mathop{\prod}\limits_{i=1}^4
(t+\tilde\eps_{i})^{-1}(t-\tilde\eps_{i}).
$$

By Proposition 1.4, $(\alpha_{1}\alpha_{2}\alpha_{3})^2\varphi$ is the 
restriction of an invariant whereas $\alpha_{1}\alpha_{2}\alpha_{3}g$ 
can be expressed via supertraces of the standard $\fs\fpe(4)$-module.

The theorem is completely proved. \qed

\begin{Corollary} Any element from $S(\fh^*)^W$ can be expressed in 
the form $\frac{P}{Q}$, where $P\in I(\fh^{*})$ and 
$Q=\mathop{\prod}\limits_{\alpha\in\tilde R_{\od}}\alpha$.
\end{Corollary}

\section*{\S 3. Anti-invariant polynomials}

The invariant polynomials appeared in \S 2 not only as elements of 
$S(\fh^*)^W$, but also as homogeneous components of supercharacters 
considered as formal power series. In this section, instead of 
supercharacters, I consider just characters and formulate (and prove) 
the corresponding analog of Chevalley's theorem.

Observe that in this case both the formulation and the proof are much 
easier.

\ssec{3.1. The two $\fg$-module structures on $U(\fg)$} It is 
well-known that the adjoint representation of $\fg$ can be uniquely 
extended to a representation in $S(\fg)$ and $U(\fg)$; moreover, the 
canonical symmetrization (actually, supersymmetrization) $\omega: 
S(\fg)\tto U(\fg)$ given by the formula
$$
\omega(x_{1}\dots x_{n})=\frac{1}{n!}
\mathop{\sum}\limits_{\sigma\in S_{n}}c(p(x), 
\sigma)x_{\sigma(1)}\dots x_{\sigma(n)},
$$
where $c(p(x), \sigma)$ is defined in (0.5.2), is an isomorphism of 
$\fg$-modules. 

Denote: $U^n(\fg)=\omega(S^n(\fg))$; we get a decomposition of $U(\fg)$
into the direct sum of $\fg$-modules.

As for the Lie algebra case, the space $U(\fg)^{*}$ can be endowed 
with the superalgebra structure by means of the coalgebra structure
($m^{*}: \fg\tto U(\fg)\otimes U(\fg)$ given by $m^{*}(x)=x\otimes 
1+1\otimes x$).

We similarly prove that $U(\fg)^{*}$ is isomorphic (as a superspace) 
to the supercommutative superalgebra of formal power series in $\dim 
\fg_{\ev}$ even indeterminates and $\dim \fg_{\od}$ odd ones.

The corresponding isomorphism is given as follows. Let $e_{1}$, 
\dots $e_{n}$ be a $\Cee$-basis of $\fg$ and $x_{1}$, \dots , $x_{n}$ 
the corresponding indeterminates of the same parity. Set:
$$
e_{\nu}=\frac{e_{1}^{\nu_{1}}}{\nu_{1}!}\dots
\frac{e_{n}^{\nu_{n}}}{\nu_{n}!},
$$
where $\nu_{1}$, \dots , $\nu_{n}$ is a collection of numbers from 
$\Zee_{+}$ for even basis vectors and 0 or 1 for odd ones. Then
$$
U(\fg)^{*}\ni f\tto S_{f}=\sum f(e_{\nu})x_{n}^{\nu_{n}}\dots
x_{1}^{\nu_{1}}.
$$
(Notice the inverse order of indeterminates.)

The {\it homogeneous component} of degree $k$ of the functional $f$ 
is a functional $f_{k}$ such that
$$
f_{k}|_{U^k(\fg)}=f|_{U^k(\fg)}\text{ and }f_{k}|_{U^l(\fg)}=0\text{ 
for }l\neq k.
$$
This splitting of $f$ into homogeneous components helps us to embed 
$S(\fg^{*})$ into $U(\fg)^{*}$ by considering $S(\fg^{*})$ as series 
whose homogeneous components of sufficiently high degree vanish.

Now, let us endow $U(\fg)$ with a $\fg$-module structure by means of 
Lemma 1.2 and $U(\fg)^{*}$ with the structure of a dual $\fg$-module.

\ssbegin{3.2}{Proposition} Consider symmetrization $\omega: 
S(\fg_{\ev})\tto U(\fg_{\ev})$ and extend it to a $\fg$-module 
homomorphism $\tilde\omega: S(\fg)\tto U(\fg)$, where $S(\fg)$ is 
considered as $\Ind_{\fg_{\ev}}^{\fg}S(\fg_{\ev})$ and $U(\fg)$ is 
considered as the $\fg$-module with respect to the structure 
$(0.5.3)$.  Then

{\em (i)} $\tilde\omega$ is a $\fg$-module isomorphism;

{\em (ii)} $\tr_{V}$ is a$\fg$-invariant element from $U(\fg)^{*}$ for any 
finite dimensional $\fg$-module $V$;

{\em (iii)} let $\tilde U^k(\fg)=\tilde\omega(\tilde S^k(\fg))$, where 
$\tilde S^k(\fg)$ is the $\fg$-submodule generated by 
$S^k(\fg_{\ev})$.  Then
$$
U(\fg)=\oplus_{k\geq 0}\tilde U^k(\fg) \eqno{(3.1)}
$$
and each of the homogeneous (with respect to (3.1)) components of an 
invariant element is an invariant itself.
\end{Proposition}

\begin{proof} (i) It suffices to show that the image of any basis of 
$S(\fg)$ under $\tilde\omega$ is a basis of $\gr U(\fg)$; but this is 
obvious.

(ii) follows form Lemma 1.2.

(iii) Follows from the fact that $\tilde U^k(\fg)$ is a $\fg$-submodule.
\end{proof}

\ssbegin{3.3}{Theorem} Let $\fg$ be one of the Lie superalgebras 
$(0.2)$, let $\fh$ be its Cartan subalgebra, $W$ its Weyl group. Then 
the restriction homomorphism onto $\fh$ induces an isomorphism
$S(\fg^{*})^{\fg}\tto S(\fh^{*})^W$.
\end{Theorem}

\begin{proof} The embedding $\fh\tto\fg$ induces the restriction 
homomorphism $U(\fg)^{*}\tto S(\fh^{*})$.  Since 
$U(\fg)=\Ind_{\fg_{\ev}}^{\fg}U(\fg_{\ev})$, there exists a bijection 
between the set of invariant elements from $U(\fg)^{*}$ and the 
$\fg_{\ev}$-invariant elements from $U(\fg_{\ev})^{*}$.

Therefore, we may assume that $\fg=\fg_{\ev}$.  If $\fg_{\ev}$ is 
reductive, the theorem is obtained as a corollary of Chevalley's 
theorem for the Lie algebras.  If $\fg_{\ev}$ is not reductive, it 
suffices to demonstrate that one can replace $\fg_{\ev}$ with 
$\fg_{0}$ for the standard $\Zee$-grading when $\fg_{0}$ is reductive 
and apply Chevalley's theorem for the Lie algebras again.
\end{proof}

\section*{Appendix.  Certain constructions with the point functor} The 
point functor is well-known in algebraic geometry since at least 1953 
\cite{W}.  The advertising of ringed spaces with nilpotents in the 
structure sheaf that followed the discovery of supersymmetries caused 
many mathematicians and physicists to realize the usefulness of the 
language of points.  F.~A.~Berezin was the first who applied the point 
functor to study Lie superalgebras.  Here we present some of his 
results and their generalizations.

All superalgebras and modules are supposed to be
finite dimensional over $\Cee $.

Thus, let $\fg$ be a Lie superalgebra, $V$ a $\fg$-module, $\Lambda $ 
the Grassmann superalgebra over $\Cee $ generated by $q$ 
indeterminates.  Define $\varphi :\Lambda \otimes V^{*}\tto \Hom 
_{\Lambda }(\Lambda \otimes V, \Lambda )$ by setting
$$
\varphi (\xi \otimes \alpha )(\eta \otimes v)=(-1)^{p(\alpha )(\eta 
)}\xi \eta \alpha (v), \; \text{for any $\xi, \eta \in \Lambda , \alpha 
\in V^{*}$.}
$$
Extend the ground field to $\Lambda $ and consider $\Lambda \otimes 
V^{*}$ and $\Hom_\Lambda(\Lambda \otimes V, \Lambda )$ as $\Lambda 
\otimes \fg$-modules.

\ssbegin{A1}{Lemma} $\varphi $ is a $\Lambda \otimes \fg$-module 
isomorphism.
\end{Lemma}

{\it Proof}.  Since $V$ is finite dimensional, $\varphi $ is a vector 
space isomorphism over $\Lambda $; besides, it is obvious that 
$\varphi $ is a $\Lambda $-module homomorphism.  Now take
$$
\xi _{1}, \xi _{2}, \xi _{3}\in \Lambda,\; \; \alpha \in V^{*},\; 
\; v\in V, \; \; x\in \fg.
$$
It is an easy exercise to prove that
$$
[(\xi _{1}\otimes x)\varphi (\xi _{2}\otimes \alpha )](\xi _{3}\otimes 
v)=\varphi (\xi _{1}\otimes x(\xi _{2}\otimes \alpha ))(\xi _{3}\otimes 
v).  \qed
$$

Consider the composition of maps
$$
V^{*}\buildrel{\varphi _{1}}\over\tto\Lambda \otimes
V^{*}\buildrel{\varphi }\over\tto\Hom _{\Lambda }(\Lambda
\otimes V, \Lambda )\buildrel{\varphi _{2}}\over\tto S_{\Lambda
}(\Hom _{\Lambda }(\Lambda \otimes V, \Lambda )), 
$$
where $\varphi _{1}(\alpha )=1\otimes \alpha$ and $\varphi _{2}$ is a 
canonical embedding of a module in its symmetric algebra.  The $\Cee 
$-module homomorphism $\varphi _{2}\circ \varphi \circ \varphi _{1}$ 
induces the algebra homomorphism
$$
S(V^{*})=S_{\Cee }(V^{*})\tto S_{\Lambda }(\Hom _{\Lambda
}(\Lambda \otimes V, \Lambda ))
$$
and, since the latter algebra is a $\Lambda $-module, we get an
algebra homomorphism
$$
\Lambda \otimes S(V^{*})\buildrel{\psi }\over\tto S_{\Lambda }(\Hom 
_{\Lambda }(\Lambda \otimes V, \Lambda )).
$$
Besides, both algebras possess a natural $\Lambda \otimes 
\fg$-module structure.

\ssbegin{A2}{Lemma}  $\psi $ {\it is a $\Lambda \otimes 
\fg$-modules and $\Lambda \otimes \fg$-algebras isomorphism}.
\end{Lemma}

\begin{proof} Let us construct the inverse homomorphism.
Consider the composition
$$
\Hom _{\Lambda }(\Lambda \otimes V, \Lambda )\buildrel{\varphi 
^{-1}}\over\tto \Lambda \otimes V^{*}\tto\Lambda
\otimes S(V^{*}).
$$
Since this composition is a $\Lambda $-module homomorphism, it
induces
the homomorphism
$$
\tilde{\psi}: S_{\Lambda }(\Hom _{\Lambda } (\Lambda \otimes V, 
\Lambda ))\tto\Lambda \otimes S(V^{*}).
$$
It is not difficult to verify that
$$
\psi \circ \tilde{\psi }|_{\Hom _{\Lambda }(\Lambda \otimes V, \Lambda 
)}=\id;\qquad \tilde{\psi }\circ \psi |_{\Lambda \otimes 
S(V^{*})}=\id;
$$
hence, $\psi$ is an isomorphism and $\tilde{\psi }$ is its inverse.  
The following proposition shows that $\psi $ is a $\Lambda \otimes 
\fg$-module isomorphism and completes the proof of Lemma A2.
\end{proof}

\ssbegin{A3}{Proposition} Let $A$, $B$ be $\Lambda $-superalgebras, 
$\fg$ a Lie superalgebra over $\Lambda $ acting by differentiations on 
$A$ and $B$.  Let $M\subset A$, $N\subset B$ be $\Lambda $-submodules 
which are at the same time $\fg$-modules generating $A$ and $B$, 
respectively, $f: A\tto B$ an algebra homomorphism such that 
$f(M)\subset N$ and $f|_{M}$ is a $\fg$-module homomorphism.  Then $f$ 
is a $\fg$-module homomorphism.
\end{Proposition}

\begin{proof} Let $a\in A$.  We may assume that $a=a_{1}\dots a_{n}$, 
where the $a_{i}\in M$.  Then for $x\in \fg$ we have
$$
\renewcommand{\arraystretch}{1.2}
\begin{gathered}
f(x(a_{1}\dots a_{n}))=f(\sum \pm a_{1}\dots xa_{i}\dots a_{n})=
\sum \pm f(a_{1})\dots f(xa_{i})\dots f(a_{n})\\
=\sum \pm f(a_{1})\dots 
xf(a_{i})\dots f(a_{n})= x[f(a_{1})\dots 
f(a_{n})]=xf(a_{1}\dots a_{n}).
\end{gathered}
$$

This proves Proposition and completes the proof of 
Lemma A2. \end{proof}

Now, let $\fh$ be a Lie superalgebra over $\Lambda $ and $U$ be a 
$\Lambda $ and $\fh$-module.  Consider $U_{\ev}$ as a $\Cee $-module.  
Then, clearly, the natural embedding $U_{\ev}\tto U$ is extendable to 
a $\Lambda $-module homomorphism $\varphi : \Lambda \otimes 
U_{\ev}\tto U$.

\ssbegin{A4}{Lemma} The homomorphism $\varphi $ is an 
$\fh_{\ev}$-module homomorphism.
\end{Lemma}

\begin{proof} Let $x\in \fh_{\ev}, \xi \in \Lambda $ and $u\in 
U_{\ev}$. Then
$$
\renewcommand{\arraystretch}{1.2}
\begin{gathered}
\varphi (x(\xi \otimes u))=\varphi (\xi \otimes xu)=\xi xu, \\
x\varphi (\xi \otimes u)=x\xi u=\text{ ( by definition of a module over a 
superalgebra) } \xi xu. \qed
\end{gathered}
$$

Thus, the adjoint map
$$
\Hom_{\Lambda }(U, \Lambda )\tto\Hom_{\Lambda }(\Lambda
\otimes U_{\ev}, \Lambda )
$$
is also an $\fh_{\ev}$-module homomorphism, therefore, by Proposition 
A3 the algebra homomorphism
$$
S_{\Lambda}(\Hom _{\Lambda}(U, \Lambda ))\tto 
S_{\Lambda}(\Hom_{\Lambda }(\Lambda \otimes U_{\ev}, \Lambda ))
$$
induced by this map is at the same time a $\fh_{\ev}$-module morphism.  
Besides, by Lemma A2 the algebra $S_{\Lambda }(\Hom _{\Lambda 
}(\Lambda \otimes U_{\ev}, \Lambda))$ is isomorphic as a $\Lambda 
\otimes \fh_{\ev}$-module and as an algebra to $\Lambda \otimes 
S(U^{*}_{\ev})$.  In particular, they are isomorphic as 
$\fh_{\ev}$-modules.

Denote by $\theta$ the composition of the homomorphisms
$$
S(V^{*})\tto\Lambda \otimes S(V^{*})\tto S_{\Lambda}(\Hom_{\Lambda 
}(\Lambda \otimes V, \Lambda ))\tto S_{\Lambda} 
(\Hom_{\Lambda}(\Lambda \otimes U_{\ev}, \Lambda )),
$$
where $U_{\ev}=(\Lambda \otimes V)_{\ev}=V_{\Lambda }$.
\end{proof} 

\ssbegin{A5}{Proposition} If $q>$dimV$_{\od}$ and $\xi \in \Lambda , 
p(\xi )=\od$, then the restriction of $\theta $ onto $\Cee [\xi 
]\otimes S(V^{*})$ is injective.
\end{Proposition}

\begin{proof} If $u\in V_{\Lambda}$, then there is defined a linear form
$L_{u}:\Hom _{\Lambda }(\Lambda \otimes V_{\Lambda }, \Lambda 
)\tto\Lambda $ by the formulas $L_{u}(l)=l(1\otimes u)$ and
$L_{u}(\xi l)=\xi l(1\otimes u)=\xi L_{u}(l)$.

Therefore, $L_{u}$ is a $\Lambda $-module homomorphism, hence, it is 
uniquely extendable to a homomorphism
$$
\varphi _{u}: \fa=S_{\Lambda }(\Hom _{\Lambda }(\Lambda \otimes 
V_{\Lambda }, \Lambda ))\tto\Lambda
$$
Consider the elements of $\fa$ as functions on $V_{\Lambda }$ setting 
$f(u)=\varphi _{u}(f)$ for $f\in \fa$ and $u\in V_{\Lambda }$.  If 
$f\in \Lambda \otimes S(V^{*})$, then set $f(u)=\varphi _{u}\circ 
\theta (f)$.  For $\alpha \in V^{*}, \xi \in \Lambda $ we have
$$
(\xi \otimes \alpha )(u)=\varphi _{u}\circ \theta (\xi \otimes \alpha 
)=L_{u}\circ \theta (\xi \otimes \alpha )=\theta (\xi \otimes \alpha 
)(1\otimes u).
$$
If $\{e_{i}\}_{i\in I}$ is a basis in $V$ and $u=\sum \lambda _{i}\otimes 
e_{i}$, 
then
$$
(\xi \otimes \alpha )(u)=\sum (-1)^{p(\alpha )p(e_{i})}\xi \lambda 
_{i}\alpha (e_{i}).  \eqno{(1)}
$$
On the other hand, the algebra $\Cee [\xi ]\otimes S(V^{*})$ is 
identified with the free supercommutative superalgebra generated by 
the $e^{*}_{i}$ and $\xi $.

Let
us assume that $p(e^{*}_{i})=0$ for $i\le n$ and $p(e^{*}_{i})=1$ for $i>n$. If
$f\in \Cee [\xi ]\otimes S(V^{*})$, then
$$
f=f_{0}+\xi f_{1}, f_{j}=\sum f_{ji_{1}\dots i_{k}}e^{*}_{i_{1}}\dots 
e^{*}_{i_{k}}, \; \text{where}\; j=0, 1\; \text{and}\; f_{ji_{1}\dots 
i_{k}}\in S(V^{*}_{\ev}).
$$
By $(1)$ we have 
$$
\renewcommand{\arraystretch}{1.2}
\begin{gathered}
f(u)=f_{0}(u)+\xi f_{1}(u)=\\
\sum f_{0i_{1}\dots i_{k}}(u)e^{*}_{i_{1}}(u)\dots 
e^{*}_{i_{k}}(u)+\sum f_{1i_{1}\dots i_{k}}(u)e^{*}_{i_{1}}(u)\dots 
e^{*}_{i_{k}}(u).
\end{gathered}
$$

Set $\lambda _{i}=a_{i}$ for $i\le n$ and $\lambda _{i}=\xi _{i-n}$ 
for $i>n$.  Then since $q>\dim V_{\od}$, we may assume that the family 
$\{\xi _{i}\}_{i\in I}$ freely generates $S(V^{*}_{\od})$ and
$$
f(u)=\sum (-1)^{k}f_{i_{1}\dots i_{k}}(a_{1}\dots a_{n})\xi 
_{i_{1}-n\dots i_{k}-n}. \eqno{(2)}
$$

If $\theta (f)=0$, then $f(u)=\varphi _{u}\circ \theta (f)$ for any 
$u\in V_{\Lambda }$.  It follows from $(2)$ that $f_{i_{1}\dots 
i_{k}}(a)=0$ for any $a\in \Cee ^{n}$.  But since $\Cee $ is 
algebraically closed, it follows (with Prop.  5.3.1 from \cite{Bu3}) 
that $f_{i_{1}\dots i_{k}}=0$; hence, $f=0$.
\end{proof} 

\ssbegin{A6}{Lemma} Let $q>$dimV$_{\od}$.  Then $f\in S(V^{*})$ is 
a $\fg$-invariant if and only if $\theta (f)\in \Lambda \otimes 
S(V^{*}_{\Lambda })$ is $\fg_{\Lambda }$-invariant.
\end{Lemma}

\begin{proof} Consider the factorization of $\theta $:
$$
\begin{gathered}
S(V^{*})\buildrel{i_{1}}\over\tto\Lambda \otimes
S(V^{*})\buildrel{i_{2}}\over\tto\\
S_{\Lambda }(\Hom _{\Lambda
}(\Lambda \otimes V, \Lambda ))\buildrel{i_{3}}\over\tto S_{\Lambda
}(\Hom _{\Lambda }(\Lambda \otimes V_{\Lambda }, \Lambda 
))\buildrel{i_{4}}\over\tto\Lambda \otimes S(V^{*}_{\Lambda }).
\end{gathered}
$$
Let $f\in S(V^{*})^\fg$, then 
$$
(\xi \otimes x)(i_{1}(f)=(\xi \otimes x)(1\otimes f)=\xi \otimes 
xf=0\; \; \text{ for }\; \xi \in \Lambda , x\in \fg.
$$

Conversely, let $yi_{1}(f)=0$ for any $y\in \fg_{\Lambda }$. 
Then
$$
0=(\xi \otimes x)(1\otimes f)=\xi \otimes xf. 
$$
If $p(y)=\od$ then let $p(\xi )=\od$. Therefore, $\xi \otimes xf=0$ and
$yf=0$.

Thus, conditions 
$$
f\in S(V^{*})^\fg\longleftrightarrow i_{1}(f)\in (\Lambda 
\otimes S(V^{*}))^\fg_{\Lambda }.
$$
Since $i_{2}$, $i_{3}$, $i_{4}$ are $\fg_{\Lambda }$-module 
homomorphisms, the above implies that if $f$ is a $\fg$-invariant then 
$\theta (f)=i_{4}\circ i_{3}\circ i_{2}\circ i_{1}(f)$ is also a 
$\fg_{\Lambda }$-invariant.

Conversely, let $\theta (f)$ be a $\fg_{\Lambda }$-invariant.  Let 
$x\in\fg_{\ev}$.  Then
$$
\theta (1\otimes xf)=\theta ((1\otimes x)(1\otimes f))=(1\otimes 
x)\theta (f)=0.
$$
 
By Proposition A5 $1\otimes xf=0$ and $xf=0$. Let $x\in 
\fg_{\od}, \; \xi \in \Lambda , \; p(\xi )=\od$. Then 
$\theta (\xi \otimes xf)= (\xi \otimes x)\theta (1\otimes f)=0$ and 
again by Proposition A5 $\xi \otimes xf=0$; hence $xf=0$ and 
therefore, $f\in S(V^{*})^\fg{{ g}}$. 
\end{proof} 

\ssec{A7. Remark} The point of the above lemmas and
propositions is that while seeking invariant polynomials on $V$
we may consider them as functions on $V_{\Lambda }$ invariant with
respect to the Lie algebra $\fg_{\Lambda }$. It makes it possible to 
apply the theory of usual Lie groups and Lie algebras and their
representations.

\ssec{A8. Remark} Let $\varphi $ be an elementary automorphism 
(of the form $\theta _{\beta }$ in Lemma 1.2.3 below) of the Lie 
algebra $\fg_{\ev}$. Clearly, $\varphi $ can be uniquely 
extended to an automorphism of the Lie superalgebra $\fg$. Let 
$\varphi (\fh)=\fh$, where $\fh$ is a Cartan subalgebra of 
$\fg$. If $i: S(\fg^{*})^{\fg}\tto S(\fh^{*})$ is the restriction 
homomorphism, then, clearly, $i(S(\fg^{*})^\fg)\subset 
S(\fh^{*})^{\varphi }$, where $A^{\varphi }$ is the set of $\varphi 
$-invariant elements of $A$.

\ssbegin{A9}{Proposition} Let $A$ be a commutative finitely 
generated algebra over $\Cee $ without nilpotents, $\fa=A\otimes 
\Lambda (p)$. Let $q\ge p$ and $f\in \fa$ be such that 
$\varphi (f)=0$ for any $\varphi :\fa\tto\Lambda (q)$. Then 
$f=0$.
\end{Proposition}

\begin{proof} Let $\psi: A\tto\Cee $ be an arbitrary homomorphism.  
Let us extend $\psi $ to a homomorphism $\varphi : \fa\tto\Lambda (q)$ 
setting $\varphi =\psi \otimes 1$.  If $\xi _{1}, \dots , \xi _{p}$ 
are generators of $\Lambda (p), f\in \fa$ and $f=\sum f_{i_{1}\dots 
i_{k}}\xi _{i_{1}}\dots \xi _{i_{k}}$, then the condition $\varphi 
(f)=0$ yields $\psi (f_{i_{1}\dots i_{k}})=0$ and, since $\psi $ is 
arbitrary, then Proposition 5.3.1 in \cite{Bu3}  shows that 
$f_{i_{1}\dots i_{k}}=0$; hence, $f=0$.
\end{proof}

\end{document}